\documentclass{elsart}

\usepackage{epsfig,psfrag}
\usepackage{amsmath}
\usepackage{amssymb}
\usepackage{latexsym,pifont,color,comment}

\newcommand{\field}[1]{\mathbb{#1}}

\def\apdq{{{\mathcal A}}^{+}\Delta{Q}}
\def\amdq{{{\mathcal A}}^{-}\Delta{Q}}
\def\Om{{\mathcal O}}
\def\reals{\field{R}}
\def\id{\mathbb I}
\def\Tm{{\mathcal T}}

\newtheorem{theorem}{Theorem}[section]

\newenvironment{proof}[1][Proof.]{\begin{trivlist}
\item[\hskip \labelsep {\bfseries #1}]}{\end{trivlist}}

\renewcommand{\qed}{\nobreak \ifvmode \relax \else
      \ifdim\lastskip<1.5em \hskip-\lastskip
      \hskip1.5em plus0em minus0.5em \fi \nobreak
      \vrule height0.75em width0.5em depth0.25em\fi}

\begin{document}

\begin{frontmatter} 

\title{A Class of Residual Distribution Schemes
and Their Relation to Relaxation Systems}

\author{James A. Rossmanith}

\address{Department of Mathematics, University of Wisconsin,
480 Lincoln Drive, Madison, WI  53706-1388, USA}

\ead{rossmani@math.wisc.edu}

\begin{abstract}
Residual distributions (RD) schemes are a class of of high-resolution
finite volume methods for unstructured grids. A key
feature of these schemes is that they make use of
genuinely multidimensional (approximate)
Riemann solvers as opposed
to the piecemeal 1D Riemann solvers usually employed by
finite volume methods. In 1D,
 LeVeque and Pelanti [{\it J. Comp. Phys.}
{\bf 172}, 572 (2001)]
 showed that
many of the standard approximate Riemann solver
methods (e.g., the Roe solver, HLL, Lax-Friedrichs)
can be obtained from applying an exact Riemann 
solver to relaxation systems of the type 
introduced by Jin and Xin [{\it Comm. Pure Appl. Math.}
{\bf 48}, 235 (1995)]. 
In this work
we extend LeVeque and Pelanti's results and obtain a multidimensional relaxation system from which multidimensional approximate
Riemann solvers can be obtained. In particular, we show that
with one choice of parameters the relaxation system yields
the standard N-scheme. With another choice,
the relaxation system yields a new Riemann solver,
which can be viewed as a genuinely multidimensional extension
of the local Lax-Friedrichs scheme. This new Riemann solver
does not require the use Roe-Struijs-Deconinck averages,
nor does it require the inversion of an $m \times m$ matrix
in each computational grid cell, where $m$ is the number of
conserved variables. Once this new scheme
is established, we apply it on a few standard cases
for the 2D compressible Euler equations of gas dynamics.
We show that
through the use of linear-preserving  limiters, the new approach produces
numerical solutions that are comparable in accuracy
to the N-scheme, despite being computationally less expensive.
\end{abstract}

\end{frontmatter}

\section{Introduction}
In the last few decades intense research into shock-capturing schemes has resulted in several numerical methods for solving partial differential equations (PDEs) that
admit discontinuous weak solutions in the form of shock-waves.
Examples of such schemes include  WENO (weighted essentially
non-oscillatory) \cite{article:JiShu96}, central \cite{article:KuTa00}, MUSCL (monotone upstream-centered
schemes for conservation laws) \cite{article:Vl79}, and wave propagation
schemes \cite{article:Le97}. One difficulty with these methods is that
in general they do not trivially extend to problems in complex
geometries. In order to handle application problems
where complex geometry is of great importance, three
broad classes of strategies have been considered:
(1) Cartesian cut-cell methods \cite{article:HeBeLe05},
(2) overlapping meshes \cite{article:ChHe90},
 and (3) unstructured grid methods.
We will focus in this work on the last approach,
thus eliminating the problem of cut-cells (1st approach) and interpolation
between different grid patches (2nd approach), but requiring some
efficient grid generation tool.

On unstructured grids, the two main classes of methods that have been
developed are discontinuous Galerkin (DG) 
\cite{article:CoShu5,article:ReedHill73} and residual distribution (RD) \cite{article:Ab01,article:De93} schemes.  The discontinuous Galerkin approach is based
on defining a piecewise polynomial approximation that is
continuous inside element interiors, but discontinuous across
element boundaries. Local 1D  Riemann problems are solved across element boundaries to construct the necessary numerical fluxes. 
Residual distribution schemes  can be viewed as a finite volume method where the finite volumes are defined by a grid that is dual
to the original triangulation. 

Although DG schemes have their own particular
advantages, the focus of this work will be on RD schemes and,
in particular, the aspect of RD schemes that separates them from
all other methods: RD schemes are based on solving genuinely multi-dimensional Riemann problems. This aspect allows one to obtain
methods that are positivity preserving for scalar  conservation laws
and essentially non-oscillatory for systems. This same feature, however,
presents a challenge: how can these multi-dimensional Riemann problems
be solved efficiently? The standard answer to this question is the so-called systems N-scheme \cite{article:Weide96} (see also
\cite{article:Ab01,article:Ab06b}), which is a generalization of Roe's approximate Riemann
solver for 1D systems \cite{article:Roe81}. One goal of this work is to
develop an alternative to this approach.

LeVeque and Pelanti \cite{article:LePe01} showed how several of the
standard approximate Riemann solvers can be interpreted as 
exact Riemann solvers for a perturbed system of hyperbolic equations
known as {\it relaxation systems}. Their work was motivated by Jin and Xin's
earlier paper \cite{article:JinXin95} on a class of numerical methods known as {\it relaxation schemes}. What LeVeque and Pelanti  essentially
showed is that Jin and Xin's ``new'' class of methods could actually be
thought of as a reinterpretation of various pre-existing approximate Riemann
solvers; these results are reviewed in Section \ref{sec:review_1d}. 
After reviewing RD schemes in Section \ref{sec:rd_schemes}, we focus in this work on the continuation of LeVeque and Pelanti's reasoning and
show how the N-scheme can be also be derived from a relaxation system.
Furthermore, using this interpretation we derive a novel genuinely 
multi-dimensional Riemann solver that can be viewed as a multidimensional extension of the 1D local Lax-Friedrichs scheme
\cite{article:Ru61}. Both of these results are presented
in Section \ref{sec:multid}. Finally, we compare the numerical accuracy
of the N-scheme and the newly derived scheme on several
examples in Section \ref{sec:num_examples}.  What we find is that
when the appropriate limiters are applied, the novel scheme has comparable accuracy to the N-scheme, although it tends to be
slightly more diffusive -- this result is of course consistent with
well-known 1D results comparing local Lax-Friedrichs versus
Roe-type approximate Riemann solvers. On the other hand, this 
loss of accuracy is compensated by the fact that the new scheme
is less computational expensive. This gain in computational efficiency
will become significant for problems involving complicated equations such as the relativistic Euler or MHD equations.
\section{Review of 1D relaxation systems}
\label{sec:review_1d}
We briefly review in this section the results
of LeVeque and Pelanti \cite{article:LePe01} for
the case of 1D conservation laws. For simplicity we
consider for the moment a scalar conservation laws
of the form
\begin{equation}
  \label{eqn:conslaw1d}
  q_{,t} + f(q)_{,x} = 0,
\end{equation}
where $x \in \reals$ is the spatial coordinate, $t \in \reals^{+}$
is the time coordinate, $q \in \reals$ is the conserved variable,
and $f(q): \reals \rightarrow \reals$ is the flux function.
We assume that this conservation law is hyperbolic,
meaning that $f'(q) \in \reals$ for all $q$ in the
solution domain.

\subsection{Finite volume methods in 1D}
Using the idea of relaxation, we will construct
in this section numerical methods for approximating
(\ref{eqn:conslaw1d}). All of these methods are
in the general class of finite volume methods
\cite{book:Le02}, which we briefly recall in this 
subsection. 

Let $\Tm_{\Delta x}$ be the numerical grid
with grid cells centered at $x = x_i$
and spanning the interval $[x_{i}-\Delta x/2, x_{i}+\Delta x/2]$,
where 
\begin{equation}
  x_i = a  +   \left( i-1/2 \right)  \Delta x.
\end{equation}
Here $i$ is an integer ranging from $1$ to $N$, $a$ and $b$ are the left and right end points of the domain, respectively, and $\Delta x = (b-a)/N$ is the grid spacing.
In each grid cell $x_i$ and at each time level $t=t^n$ we seek 
an approximation to the cell average of the exact solution
$q(x,t)$:
\begin{equation}
  Q^{n}_{i} \approx \frac{1}{\Delta x} \int_{x_i - \Delta x/2}^{x_i + \Delta x/2}
    q(\xi,t^n) \, d\xi.
\end{equation}
Integrating (\ref{eqn:conslaw1d}) over the
grid cell centered at $x_i$ and from $t=t^n$ to $t=t^{n+1}$
results in a numerical update formula for $Q^n_i$ that
can be written in the following {\it fluctuation splitting form}:
\begin{equation}
   \label{eqn:update1d}
  Q^{n+1}_i = Q^n_i - \frac{\Delta t}{\Delta x}
  \left[ \amdq^n_{i+1/2} + \apdq^n_{i-1/2} \right],
\end{equation}
where $\amdq^n_{i+1/2}$ and $\apdq^n_{i-1/2}$ 
are left- and right-going {\it fluctuations}, which measure the amount of 
flux that enters into grid cell $x_i$ through the grid interfaces
at $x=x_i + \Delta x/2$ and $x=x_i  - \Delta x/2$, respectively.
In order for this update to be numerically conservative
these fluctuations must satisfy
\begin{equation}
  \label{eqn:cons1d}
  \amdq^n_{i+1/2} + \apdq^n_{i+1/2} = f(Q_{i+1}^n) 
  - f ( Q_{i}^{n} ).
\end{equation}
Note that in update (\ref{eqn:update1d}) we collect
the left-going fluctuation from the grid interface 
at $x_{i}+\Delta x/2$ and the right-going fluctuation from 
the grid interface at $x_{i}-\Delta x/2$, while in 
expression (\ref{eqn:cons1d}) we are adding the left- and right-going
fluctuations at the same grid interface.

For first-order accurate methods, the fluctuations in update (\ref{eqn:update1d}) are obtained by first assuming that 
the approximate solution has a constant value, $Q_i^{n}$,
in each grid cell, and then solving at each grid interface, 
$x_{i-1/2} \equiv x_{i} - \Delta x/2$, the initial value problem for (\ref{eqn:conslaw1d})
with the piecewise constant initial data:
\begin{equation}
    q(x,0) = \begin{cases}
    		Q_{i-1}^n & \text{if }  x < x_{i-1/2}, \\
		Q_{i}^n & \text{if }  x < x_{i-1/2}.
    	\end{cases}
\end{equation}
This initial value problem is referred to as the
{\it Riemann problem}. One of the pieces
of information that can be obtained from
solving the Riemann problem is how
much of the initial flux difference, $f(Q_{i}^n) - f(Q_{i-1}^n)$,
is carried to the left and how much to the right. It is precisely
this information that is stored in the fluctuations,
$\amdq$ and $\apdq$.

\subsection{Relaxation method framework in 1D}
The relaxation schemes introduced by Jin
and Xin \cite{article:JinXin95} are based on the idea of
approximating the quasilinear equation (\ref{eqn:conslaw1d}) by a 
linear system with a cleverly chosen source term. The role
of this source term is to force the linear system to {\it relax} in the limit as $t \rightarrow \infty$ towards the original equation.
By ``hiding'' the nonlinearity in the source term, relatively
complicated quasilinear Riemann problems can be
replaced by simpler linear Riemann problems. 

There are many kinds of relaxation systems
that one could develop in order to create
an approximate solution to (\ref{eqn:conslaw1d})
(see pp. 26--48 of 
Bouchut \cite{book:Bouchut05} for a discussion of several
different approaches). In this work we follow
the approach of \cite{article:LePe01} and
consider the following relaxation system:
\begin{gather}
\label{eqn:relax_scalar_1d}
\begin{bmatrix}
  q \\ \mu 
\end{bmatrix}_{,t} + 
\underset{\text{coefficient matrix}}{\underbrace{
\begin{bmatrix}
  0 \quad &  1 \\
  - c \, d  \quad &  c + d
\end{bmatrix}}}
\begin{bmatrix}
  q \\ \mu
\end{bmatrix}_{,x} = \frac{1}{\varepsilon} 
\underset{\text{source term}}{\underbrace{\begin{bmatrix}
 0 \\
 f(q) - \mu 
\end{bmatrix}}},
\end{gather}
where $c, d \in \reals$ are parameters that will
be adjusted in the next few subsections in order to arrive
at various approximate Riemann solvers. Without loss of 
generality we will assume that $c \le d$.
The key observation is that by taking $\varepsilon \rightarrow 0$,
the right-hand side forces $\mu \rightarrow f(q)$. Since
the first equation in the above system is $q_{,t} + \mu_{,x} = 0$,
$\mu \rightarrow f(q)$ will cause the relaxed system solution to approach
the original conservation law solution.

In order to make this statement more precise, we will
carry out a so-called {\it Chapman-Enskog} expansion,
which in this case is simply a Taylor series expansion in
$\varepsilon$ applied to system (\ref{eqn:relax_scalar_1d}).
Omitting the algebra, this expansion 
to $\Om(\varepsilon^2)$ yields the following equation for
$q(x,t)$:
\begin{equation}
  \underset{\text{original cons. law}}{\underbrace{q_{,t} + f(q)_{,x}}} = \underset{\text{diffusive correction}}{\underbrace{\varepsilon \left[ \left( \frac{\partial f}{\partial q}
   - c \right) \left( d - \frac{\partial f}{\partial q} \right) q_{,x}
   \right]_{,x}}} + {\mathcal O}(\varepsilon^2).
\end{equation}
This approximation is stable for $\varepsilon >0$ if
the values $c$ and $d$ are chosen to produce
positive (or at least non-negative) diffusion; this occurs if
\begin{equation}
   c \le \frac{\partial f}{\partial q} \le d.
\end{equation}
The above statement is often referred
as the {\it sub-characteristic condition}
(see for example \cite{article:ChLeLi94,article:JinXin95,article:LePe01}),
since it requires that  the eigenvalues of the
coefficient matrix, which are just $c$ and
$d$, enclose the characteristic speed of
the original conservation law, $\partial f/\partial q$.

From relaxation system (\ref{eqn:relax_scalar_1d}), LeVeque
and Pelanti \cite{article:LePe01} showed that 
various classical approximate Riemann solvers could
be derived. Following the philosophy of operator splitting
 (see pp. 380---390 of LeVeque \cite{book:Le02}
 for a review), system (\ref{eqn:relax_scalar_1d})
 is first rewritten as two sub-problems:
 \begin{gather}
\label{eqn:relax_scalar_1d_a}
\begin{bmatrix}
  q \\ \mu 
\end{bmatrix}_{,t} + 
\begin{bmatrix}
  0 \quad &  1 \\
  -c \, d  \quad &  c + d
\end{bmatrix}
\begin{bmatrix}
  q \\ \mu
\end{bmatrix}_{,x} = 0, \\
\label{eqn:relax_scalar_1d_b}
\mu_{,t} = \frac{1}{\varepsilon} \left( f(q) - \mu \right).
\end{gather}
Using this interpretation, LeVeque and Pelanti's \cite{article:LePe01} procedure for obtaining different approximate Riemann solvers can be summarized as follows:
\begin{enumerate}
  \item Choose values for the parameters $c$
  	and $d$.
  \item Exactly solve the Riemann problem for the
	   homogeneous linear system 
	   (\ref{eqn:relax_scalar_1d_a}).
  \item Approximate the effect of equation
  	 (\ref{eqn:relax_scalar_1d_b})
  	on the solution calculated in Step (2)
  	by directly setting  $\mu = f(q)$. In other words,	
 	instantaneously relax the solution from Step (2)
	to the $\varepsilon \rightarrow 0$ limit.
\end{enumerate}
We will simply refer to this as the {\it relaxation procedure}.
In the next four subsections, we will apply this strategy for
various values of $c$ and $d$. Each time we carry out
step (2) of the above procedure we will exactly solve
the initial value problem (i.e., $-\infty < x < \infty$) 
for system (\ref{eqn:relax_scalar_1d_a}) using
the generic {\it Riemann data}:
\begin{equation}
  q(x,0) = \begin{cases}
  	Q_{\ell} & \text{if }  x < 0, \\
	Q_{r} & \text{if }  x > 0,
	\end{cases}
	\quad \text{and} \quad
 \mu(x,0) = \begin{cases}
  	f(Q_{\ell}) & \text{if }  x < 0, \\
	f(Q_{r})& \text{if }  x > 0,
	\end{cases}
\end{equation}
where $Q_{\ell}$ and $Q_{r}$ are constants.
Note that we are allowed to take $\mu = f(q)$
in the initial conditions at any arbitrary time step, since in the
previous time-step we set $\mu = f(q)$ in Step (3)
of the relaxation procedure.

\subsection{Local Lax-Friedrichs (LLF) for scalar equations}
\label{sec:LLF1D}
The local Lax-Friedrichs or Rusanov method 
\cite{article:Ru61} is obtained by applying the
relaxation procedure with the choice
\begin{equation}
  c = -s \quad \text{and} \quad d = s,
\end{equation}
where $s \ge | f'(q) |$ in order to satisfy the sub-characteristic
condition.
With this choice the Riemann solution is
obtained by splitting the jump between the left and right states,
$(Q_{\ell}, f(Q_{\ell}))$ and $(Q_r, f(Q_r))$, along
the eigenvectors of the coefficient matrix:
\begin{equation}
\begin{bmatrix}
	Q_r - Q_{\ell} \\
	f(Q_r) - f(Q_{\ell})
\end{bmatrix}
 = \alpha^1 \begin{bmatrix} 1 \\ -s \end{bmatrix}
 +  \alpha^2 \begin{bmatrix} 1 \\ s \end{bmatrix},
\end{equation}
where the corresponding eigenvalues are $\lambda^1 = -s$
and $\lambda^2 = s$.
From this expression we obtain the following fluctuations:
\begin{align}	
  \label{eqn:LLF1}
  \amdq &= \lambda^1 \, \alpha^1 = 
  	\frac{1}{2} \left( f(Q_r) - f(Q_{\ell}) \right) - \frac{s}{2} \left( Q_r - Q_{\ell} \right), \\
	\label{eqn:LLF2}
  \apdq &=  \lambda^2 \, \alpha^2 =\frac{1}{2} \left( f(Q_r) - f(Q_{\ell}) \right) + \frac{s}{2} \left( Q_r - Q_{\ell} \right).
\end{align}

\subsection{Harten, Lax, and van Leer (HLL) for scalar equations}
The HLL method of
\cite{article:HaLaLe83} is obtained by applying the
above procedure with the choice
\begin{equation}
  c = s_{\ell} \quad \text{and} \quad d = s_r,
\end{equation}
where $s_{\ell} \le f'(q) \le s_r$ in order to satisfy the sub-characteristic
condition.
With this choice the Riemann solution is
obtained by splitting the jump between the left and right states,
 $(Q_{\ell}, f(Q_{\ell}))$ and $(Q_r, f(Q_r))$, along
the eigenvectors of the coefficient matrix:
\begin{equation}
\begin{bmatrix}
	Q_r - Q_{\ell} \\
	f(Q_r) - f(Q_{\ell})
\end{bmatrix}
 = \alpha^1 \begin{bmatrix} 1 \\  s_{\ell}  \end{bmatrix}
 +  \alpha^2 \begin{bmatrix} 1 \\ s_r \end{bmatrix},
\end{equation}
where the corresponding eigenvalues are $\lambda^1 = s_{\ell}$
and $\lambda^2 = s_r$.
From this expression we obtain the following fluctuations:
\begin{gather}
  \begin{split}
  \label{eqn:HLL1}
  {\mathcal A}^{\pm} \Delta Q  &= 
   s_{\ell}^{\pm} \, \alpha^1 +  s_{r}^{\pm} \, \alpha^2 \\
  	&=\left( \frac{s_{r}^{\pm} - 
	s_{\ell}^{\pm} }{s_r - s_{\ell}} \right)
	 \left( f(Q_r) - f(Q_{\ell}) \right)
	 - \left( \frac{s_{r}^{\pm} s_{\ell} - 
	s_{\ell}^{\pm} s_{r}}{s_r - s_{\ell}} \right)
	\left( Q_r - Q_{\ell} \right).
	\end{split}
\end{gather}
In the above expressions we have made use of the 
following notation:
\begin{equation}
     s^{+} = \text{max}(0, s) \quad \text{and} \quad
     s^{-} = \text{min}(0, s).
\end{equation}
We will make use of this notation throughout the
remainder of this paper.

\subsection{Roe's approximate Riemann solver for scalar equations}
Roe's approximate Riemann solver
\cite{article:Roe81} is obtained by applying the
above procedure with the choice
\begin{equation}
  c = d = s \equiv \frac{f(Q_r) - f(Q_{\ell})}{Q_r - Q_{\ell}}.
\end{equation}
With this choice the coefficient matrix becomes deficient, since
only one linearly independent eigenvector exists.
Therefore, the jump between the left and right states,
$(Q_{\ell}, f(Q_{\ell}))$ and $(Q_r, f(Q_r))$, can be
written as
\begin{equation}
\begin{bmatrix}
	Q_r - Q_{\ell} \\
	f(Q_r) - f(Q_{\ell})
\end{bmatrix}
 = \alpha \begin{bmatrix} 1 \\  s  \end{bmatrix},
\end{equation}
where the corresponding eigenvalue is $\lambda = s$
(algebraic multiplicity 2, geometric multiplicity 1).
Although this seems like an over-determined system
for $\alpha$, there exists a unique solution:
$\alpha = Q_r - Q_{\ell}$. This
results in the following fluctuations:
\begin{align}
   {\mathcal A}^{\pm} \Delta Q &= s^{\pm} \left( Q_r - Q_{\ell} \right).
\end{align}

\subsection{LLF and HLL for systems}
Finally, we briefly explain how these three interpretations
can be applied to a system of conservation laws
of the form (\ref{eqn:conslaw1d}), now with  $q \in \reals^{m}$ and
$f(q): \reals^m \rightarrow \reals^m$.  We again assume
hyperbolicity, which implies that the $m \times m$ matrix,
$\partial f/\partial q$, has $m$ real eigenvalues and 
$m$ linearly independent eigenvectors for all $q$ in the solution domain.

The systems LLF and HLL methods are obtained
by considering the following relaxation system:
\begin{equation}
 \label{eqn:relax_sys1}
  \begin{bmatrix} q \\ \mu \end{bmatrix}_{,t} + 
  \begin{bmatrix} 0 \id \quad & \id \\
  	- c \, d \, \id \quad & \left( c + d \right) \id
   \end{bmatrix} \begin{bmatrix} q \\ \mu \end{bmatrix}_{,x}
   =\frac{1}{\varepsilon} 
\begin{bmatrix}
 0 \\
 f(q) - \mu
\end{bmatrix},
\end{equation}
where $\id$ is the $m \times m$ identity matrix, 
$0 \id$ is the $m \times m$ matrix with zeros in
every entry, $\mu \in \reals^m$, and $c, d \in \reals$. 

The systems LLF method is obtained by taking
\begin{equation}
  s = d = -c, \quad \text{where} \quad s \ge \max_{p=1, \ldots, m}  | \lambda^{p} |,
\end{equation}
and $\lambda^{p}$ is the $p^{\text{th}}$ eigenvalue
of $\partial f/\partial q$.
With this choice we again arrive at
formula (\ref{eqn:LLF1}), which is now
applied to each component of the solution vector.

Similarly, the systems HLL method is obtained by taking
\begin{equation}
  s_{\ell} = c, \quad
  s_r = d, \quad \text{where} \quad 
 s_{\ell}  \le \min_{p=1, \ldots, m}  \left( \lambda^{p} \right)
\quad \text{and} \quad
  s_{r}  \ge \max_{p=1, \ldots, m} \left( \lambda^{p} \right).
\end{equation}
With this choice we again arrive at
formula (\ref{eqn:HLL1}), which is now
applied to each component of the solution vector.

\subsection{Roe's approximate Riemann solver for systems}
Roe's approximate Riemann solver does not follow from
working with relaxation system (\ref{eqn:relax_sys1}), but instead from \begin{equation}
 \label{eqn:relax_sys2}
  \begin{bmatrix} q \\ \mu \end{bmatrix}_{,t} + 
  \begin{bmatrix} 0 \id \quad & \id \\
  	- (\hat{J})^2 \quad & 2 \hat{J} 
   \end{bmatrix} \begin{bmatrix} q \\ \mu \end{bmatrix}_{,x}
   =\frac{1}{\varepsilon} 
\begin{bmatrix}
 0 \\
 f(q) - \mu
\end{bmatrix}.
\end{equation}
In the above expression,
\begin{equation}
 \hat{J}  =  \frac{\partial f}{\partial q} \left( \hat{Q} \right),
\end{equation}
where $\hat{Q}$ is the Roe average \cite{article:Roe81}
and satisfies
\begin{equation}
\label{eqn:roe_ave}
\frac{\partial f}{\partial q} \left( \hat{Q} \right) \left( Q_r - Q_{\ell} \right) = f(Q_r) - f(Q_{\ell}).
\end{equation}
With this choice the Riemann solution is obtained by splitting the
jump between the left and right states, $(Q_{\ell}, f(Q_{\ell}))$ 
and $(Q_r, f(Q_r))$, along the $m$ distinct eigenvectors of  $\hat{J}$
in the following manner:
\begin{equation}
\label{eqn:roe_rp}
\begin{bmatrix}
	Q_r - Q_{\ell} \\
	f(Q_r) - f(Q_{\ell})
\end{bmatrix}
 = \alpha^1 \begin{bmatrix} r^1 \\ s^1 r^1  \end{bmatrix} 
+ \cdots + \alpha^m \begin{bmatrix} r^m \\  s^m r^m  
\end{bmatrix},
\end{equation}
where $s^p$ and $r^p$ are the $p^{\text{th}}$ eigenvalue
and right eigenvector of the Roe matrix $\hat{J}$, respectively.
Just as in the scalar case, it seems as though the parameters
$\alpha$ are overdetermined. However, since $\hat{J}$ satisfies
the constraint (\ref{eqn:roe_ave}), it can easily be shown that
the first set of $m$ equations involving
$Q_r - Q_{\ell}$ are identical to the second set of $m$ equations
involving $f(Q_r) - f(Q_{\ell})$. In other words, there are only $m$
distinct equations for $m$ values of $\alpha$; and therefore,
a unique solution exists:
\begin{equation}
   \alpha^p = \ell^p \cdot \left( Q_r - Q_{\ell} \right), \quad
  \text{for} \quad p = 1, \ldots, m,
\end{equation}
where $\ell^{p}$ is the $p^{\text{th}}$ left eigenvector of $\hat{J}$.
This results in the following fluctuations:
\begin{align}
\label{eqn:Roe_fluct1}
  {\mathcal A}^{\pm} \Delta Q  &= \sum_{p=1}^m  \left[ s^p \right]^{\pm} \,  \left\{ \ell^p \cdot \left( Q_r - Q_{\ell} \right) \right\}  r^p.
\end{align}
Note that conservation follows from (\ref{eqn:roe_ave}).

\section{Residual distribution schemes}
\label{sec:rd_schemes}
We will describe in this section the basic residual
distribution method for solving hyperbolic conservation
laws in multidimensions. For further details we refer
the reader to articles by Abgrall \cite{article:Ab01,article:Ab06} and Abgrall and Mezine \cite{article:AbMe03,article:AbMe04}.
We consider a conservation law of the form
\begin{equation}
  \label{eqn:multid_conslaw}
    q_{,t} + \nabla \cdot \vec{f}(q) =  q_{,t} + f(q)_{,x} + g(q)_{,y} = 0,
\end{equation}
where $(x,y) \in \reals^2$ are the spatial coordinates,
$t \in \reals$ is the time coordinate, $q \in \reals^m$ is
the vector of conserved variable, and $\vec{f}(q): \reals^m \rightarrow \reals^{m \times 2}$
is the flux function.
We will assume that this equation is hyperbolic, meaning
that the $m \times m$ flux Jacobian matrix,
\begin{equation}
   J(\vec{n}) \equiv \vec{n} \cdot \frac{\partial \vec{f}(q)}{\partial q},
\end{equation}
is diagonalizable with real eigenvalues for all $\vec{n} \in \reals^2$
such that $\| \vec{n} \| = 1$
and for all $q$ in the solution domain.

\begin{figure}[!e]
\begin{center}
\psfrag{Triangle center}{Triangle center}
\psfrag{Triangle node}{Triangle node}
\psfrag{Grid lines of median dual cells}{Grid lines of median dual cells}
\psfrag{Grid lines}{Grid lines}
\psfrag{Center of triangle edge}{Center of triangle edge}
\includegraphics[height=84mm]{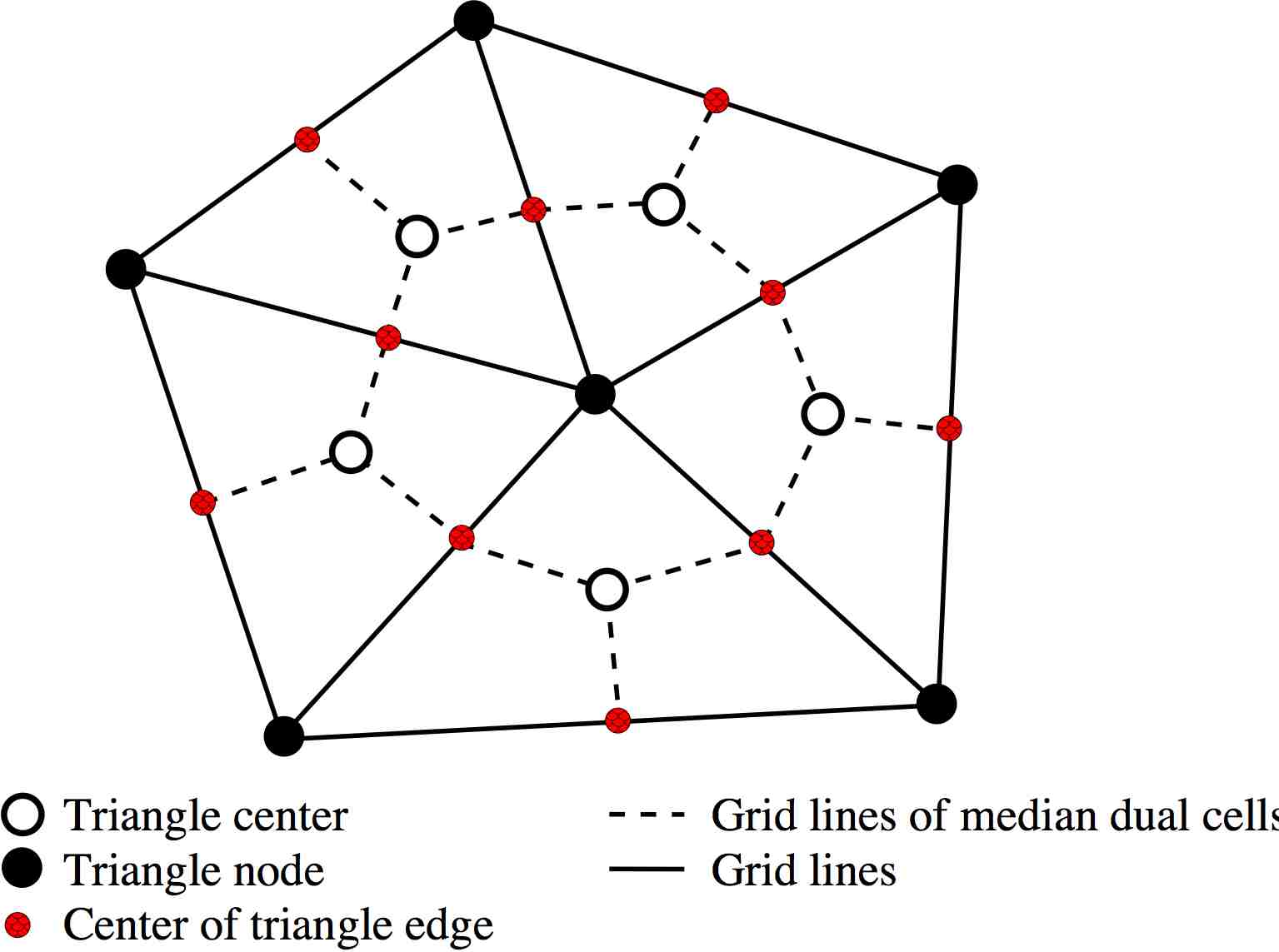}
\end{center}
\caption{Sample triangulation and dual grid.}
\label{fig:grid_dual_grid}
\end{figure}

The first step in approximately solving (\ref{eqn:multid_conslaw})
in some domain $\Omega \subset \reals^2$ is to 
mesh the domain with a finite number of triangles. 
We will refer to this triangulation as $\Tm_{h}$, where
$h$ refers to a representative triangle radius, which in this
work we just take to be the square root of the triangle area:
$h \equiv \sqrt{|\Tm|}$.
Associated with this triangulation is a {\it dual grid}, which is
constructed by connecting triangle centers to edge centers.
A an example triangulation along with its dual grid is shown in Figure
\ref{fig:grid_dual_grid}.

\begin{figure}[!t]
\begin{center}
\psfrag{Q1}{$Q_1$}
\psfrag{Q2}{$Q_2$}
\psfrag{Q3}{$Q_3$}
\psfrag{n1}{$\vec{n}_1$}
\psfrag{n2}{$\vec{n}_2$}
\psfrag{n3}{$\vec{n}_3$}
\includegraphics[height=50mm]{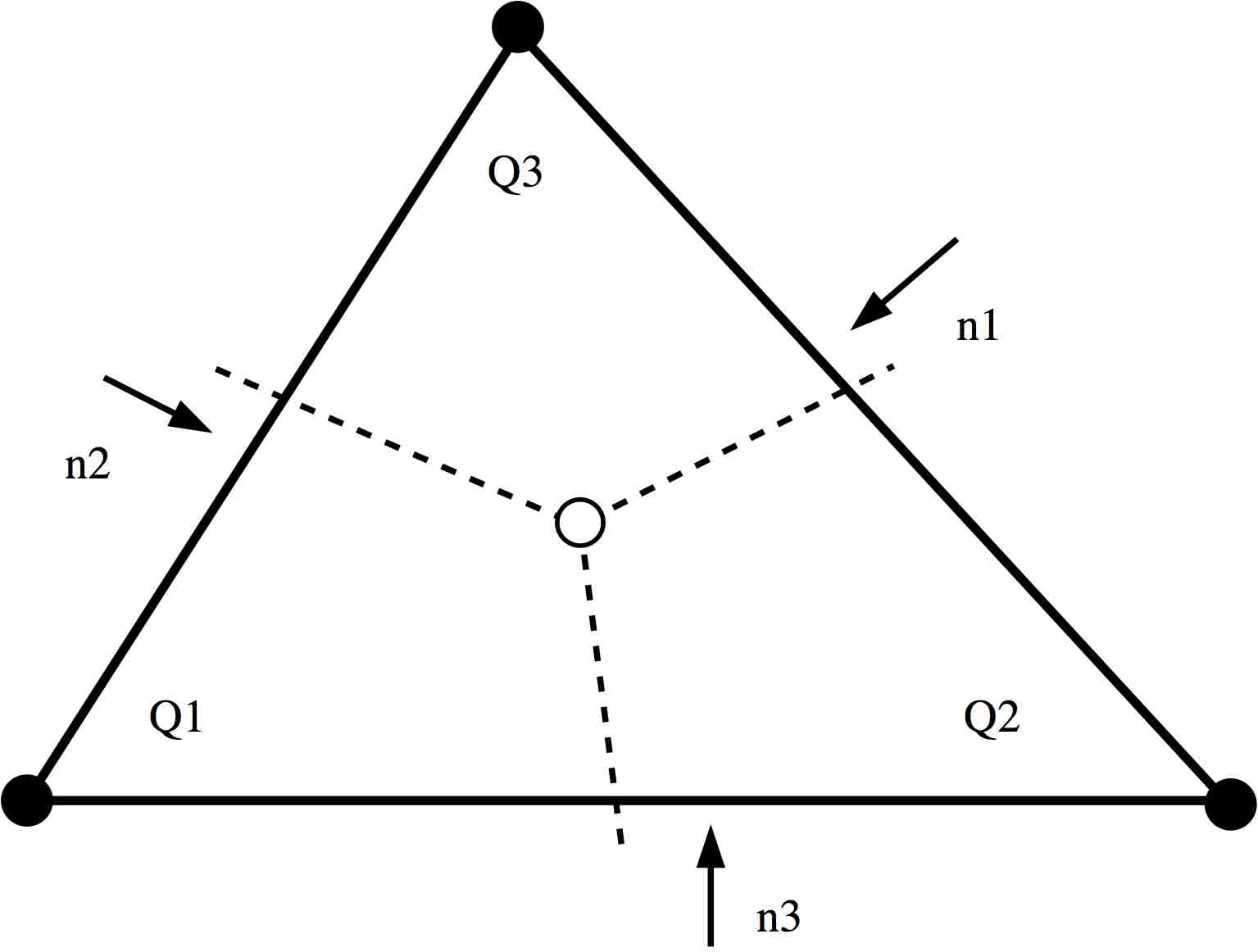}
\end{center}
\caption{A depiction of the multidimensional Riemann problem that
must be solved in each triangle. The numerical solution is piecewise
constant on each median dual cell. For example, the approximate
solution on the three dual
cells that overlap the triangle shown in this figure are
$Q_1$, $Q_2$, and $Q_3$.
Note that the area of each of the three
sections is the same, the midpoint where the dashed lines
meet is $\left( \vec{x}_1 + \vec{x}_2 + \vec{x}_3 \right)/3$,
$\vec{n}_k$ are the inward-pointing normal vectors 
to each edge, and the magnitude of $\vec{n}_k$ is equal
to the length of the edge to which it is orthogonal.
}
\label{fig:multid_riemann}
\end{figure}

Unlike the discontinuous Galerkin approach
 \cite{article:CoShu5,article:ReedHill73},
approximate solutions are
 centered on triangle nodes (i.e., centers of the dual
 grid) rather than triangle centers. 
In order to obtain an update for these node centered values,
we integrate (\ref{eqn:multid_conslaw}) over the median
dual cell $C_i$ and from $t=t^n$ to $t=t^{n+1}$:
\begin{equation*}
\begin{split}
    \iint_{C_i} q(\vec{x}, t^{n+1}) \, d\vec{x}
    &= \iint_{C_i} q(\vec{x}, t^{n}) \, d\vec{x} - \int_{t^n}^{t^{n+1}} \iint_{C_i} \nabla \cdot \vec{f}(q) \, d\vec{x} \, dt \\
    &=  \iint_{C_i} q(\vec{x}, t^{n}) \, d\vec{x}
    - \sum_{\Tm: \, i \in \Tm} \int_{t^n}^{t^{n+1}} \oint_{\partial \left( C_i \cap \Tm \right)}
    \vec{f}(q) \cdot d\vec{s} \, dt.
    \end{split}
\end{equation*}
Next we define the median dual cell average and
the time-averaged fluctuations
through $\partial \left( C_i \cap \Tm \right)$:
\begin{align}
   Q^{n}_i &\equiv \frac{1}{| C_i |}
   \iint_{C_i} q(\vec{x}, t^{n}) \, d\vec{x}, \\
   \Phi^{\Tm}_i &\equiv \frac{1}{\Delta t}
      \int_{t^n}^{t^{n+1}} \oint_{\partial \left( C_i \cap \Tm \right)}
    \vec{f}(q) \cdot d\vec{s} \, dt.
\end{align}
Using these definitions, the update formula for a 
generic residual distribution scheme can be written as follows:
\begin{equation}
\label{eqn:forwardeuler}
Q^{n+1}_i  =  Q^n_i - \frac{\Delta t}{| C_i |}
   \sum_{\Tm: \, i \in \Tm} \, \Phi^{\Tm}_i.
\end{equation}
In the remainder of this work, we will approximate the
exact solution, $q(\vec{x},t)$, with a piecewise constant
representation, $Q^{n}_i$, that is constant on each median
dual cell. We note that this view of RD schemes is slightly
different than the standard view (e.g., \cite{article:Ab01}), where
the approximate solution is usually viewed to be piecewise linear
on each triangle $\Tm$. Although these descriptions seem
contradictory, in the case for first-order accuracy in time,
both interpretations yield the same numerical schemes.
The advantage of viewing the solution as being piecewise
constant on each medial dual cell is that this naturally
sets up a series of multidimensional Riemann problems
in each triangle (see Figure \ref{fig:multid_riemann}), which
can be solved to construct the fluctuations\footnote{In this work, 
the terms {\it distributed residual} and {\it fluctuation} mean the
same thing and are used interchangeably.} $\Phi_i^{\Tm}$.
In this way, we can then view approximate constructions
of  $\Phi_i^{\Tm}$ as approximate Riemann solvers.

Computing the fluctuations $\Phi_i^{\Tm}$ is generally
done using the following framework (again, the two
interpretations, piecewise
constant on each dual cell vs. piecewise linear on each
primal cell, make use of the exactly the same framework):
\begin{enumerate}
   \item On each triangle $\Tm$ construct a
   		total residual:
		\begin{equation}
		\label{eqn:tot_res}
		  \Phi^{\Tm} = \iint_{\Tm} \nabla \cdot \vec{f}^{\, h} \, d\vec{x}
		  = \oint_{\partial \Tm} \vec{f}^{\, h} \cdot \, d\vec{s},
		\end{equation}
		where $\vec{f}^{\, h}$ denotes an interpolant
	 	that passes through the three nodal values
		\begin{equation*}
		 \left( \vec{x}_i, \, \vec{f}\left(Q_i\right) \right)
		 \quad \text{for} \quad
		 i=1,2,3.
		 \end{equation*}
		 For example if we simply use linear
		 interpolation, the total residual can be
		 written as\footnote{For the standard N-scheme, which we
		 will describe shortly, this is not the 
		 definition of the total residual that is used.}
		 \begin{equation}
		 	\label{eqn:tot_res_linear}
		    \Phi^{\Tm} = \frac{1}{2} \sum_{i=1}^{3} 
		    \vec{f}\left( Q_i \right) \cdot \vec{n}_i.
		 \end{equation}
		 Here $\vec{n}_i$ represents the inward pointing
		 normal vectors to the three edges of the triangle $\Tm$.
		 The length of $\vec{n}_i$ is equal to the length of the
  	 	edge to which it is perpendicular. If the three nodes of
		triangle $\Tm$ are given by $(x_i, y_i)$ for $i=1,2,3$, then
		the three scaled normals can be written as
		\begin{align*}
		    \vec{n}_1 &= \left(  y_2 - y_3, \, x_3 - x_2 \right)^t, \\
		    \vec{n}_2 &= \left(  y_3 - y_1, \, x_1 - x_3 \right)^t, \\
		    \vec{n}_3 &= \left(  y_1 - y_2, \, x_2 - x_1 \right)^t.
		\end{align*}
  \item Once this total residual has been calculated, it is then
  	distributed to each of the nodes of the triangle:
	\begin{equation*}
	  \Phi^{\Tm} \quad \rightarrow \quad \Phi^{\Tm}_1, \, \, \Phi^{\Tm}_2, \, \, \Phi^{\Tm}_3.
	\end{equation*}
	The detailed strategy for how this distribution is accomplished
	yields a specific numerical method.
\end{enumerate}

\subsection{Design principles for scalar conservation laws}
We first focus on design principles for scalar equations;
in a subsequent subsection we explain how to extend this
to the systems case.
In order to obtain a numerical update that
produces a stable and accurate approximation 
to (\ref{eqn:multid_conslaw}), we will need the distribution strategy
to satisfy certain properties:
\begin{enumerate}
  \item {\bf Numerical conservation:} Since the interpolation of
	the numerical solution is continuous across triangle
	edges, conservation simply reduces to the following
	constraint:
	\begin{equation}
		\label{eqn:cons_constraint}
	    \sum_{i=1}^{3} \Phi^{\Tm}_i =  \Phi^{\Tm}.
	\end{equation}
	In other words, in a given triangle, 
	the sum of the distributed residuals
	must equal the total residual.
	
  \item  {\bf Monotonicity preserving:} 
  	This condition makes sure that the numerical
	update satisfies a local maximum principle,
	which is needed to guarantee that the update
	does not generate any new spurious maxima
	or minima.
	If we write
	\begin{equation}
	\label{eqn:monotone_condition}
	  \Phi^{\Tm}_i = \sum_{j=1}^{3}  c^{\Tm}_{ij} \left( Q^n_i - Q^n_j \right),
	\end{equation}
	then the monotonicity requirement can be written as
	(see \cite{article:De93}):
	\begin{equation}
	   c^{\Tm}_{ij} \ge 0 \quad \forall \, i,j \in \left( 1, 2, 3 \right).
	\end{equation}
	
  \item {\bf Linear preserving}: The order of accuracy of
  	update (\ref{eqn:forwardeuler}) in the steady-state
	 depends, among other things, on how accurately the total
	 residual (\ref{eqn:tot_res}) is calculated on each triangle
	 \cite{article:AbRo03}. If we use formula (\ref{eqn:tot_res_linear}),
	 then $\Phi^{\Tm} = {\mathcal O}(h^3)$ in the steady state,
	 which is the correct order of accuracy if want an
	 approximate solution, $Q_i$,
	 that is ${\mathcal O}(h^2)$ accurate in the steady state.
	 What we actually need in order to get an ${\mathcal O}(h^2)$
	 accurate steady state solution is that not only that
	 $\Phi^{\Tm} = {\mathcal O}(h^3)$, but that each
	 distributed residual also satisfies $\Phi^{\Tm}_i
	 = {\mathcal O}(h^3)$.
	 The distributed residuals can be written as
	 \begin{equation}
	 	\label{eqn:weights}
	    \Phi^{\Tm}_1 = \beta_1 \Phi^{\Tm}, \quad
	    \Phi^{\Tm}_2 = \beta_2 \Phi^{\Tm}, \quad
	    \Phi^{\Tm}_3 = \beta_3 \Phi^{\Tm},
	 \end{equation}
	 where $\beta_i$ measures the fraction of the
	 total residual that is distributed to node $i$. 
	  To ensure that the distributed residuals are
	  of the same order as the total residual, 
	   we need to make sure that
	  the $\beta_i$'s remain bounded as $h \rightarrow 0$.
	  Therefore,
	  the ${\mathcal O}(h^2)$ accuracy condition, or
	  more commonly referred to as the linear preserving
	  condition, can be writtten as follows (see \cite{article:De93}):
	  \begin{equation}
	    \beta_i \quad \text{for $i=1,2,3$, 
	    is uniformly bounded independent
	    of the mesh}.
	  \end{equation} 	 
\end{enumerate}

\subsection{Scalar N-scheme}
Modern finite volume methods for hyperbolic PDEs
are typically based on solving, either exactly or approximately,
a Riemann problem between neighboring states. For
multidimensional problems, a standard approach is to 
solve local 1D Riemann problems and then use the
information from the Riemann solutions to construct
numerical fluxes or fluctuations (see Chapters 19-21 of
LeVeque \cite{book:Le02}). 

In the RD framework, however,
a multidimensional Riemann problem is solved. In an arbitrary
triangle $\Tm$, we consider the Riemann problem between
three constant states: $Q_1$, $Q_2$, and $Q_3$ (see
Figure \ref{fig:multid_riemann} for a depiction). Exact
solutions to multidimensional Riemann problems are
at best expensive to evaluate, and in general not well-understood
for many hyperbolic systems such as the Euler equations from
gas dynamics \cite{article:ScCoGl93}.
Therefore, in practice an approximate method such as
the N-scheme (the ``N'' stands for Narrow) \cite{article:Decon93,article:De93} is
utilized; this approach can be viewed
as a multidimensional generalization of Roe's approximate
Riemann solver \cite{article:Roe81}.

Just as in the 1D case, we define a Roe-average
(henceforth called the Roe-Struijs-Deconinck average \cite{article:Decon93}):
\begin{equation}
\vec{u}^{\, \Tm} \equiv \frac{\partial \vec{f}}{\partial q} \left(\hat{Q}\right),
\end{equation}
where $\hat{Q}$ is an average of
the three nodal values $Q_i$ for $i=1,2,3$ on the
current triangle $\Tm$.
The Roe-Struijs-Deconinck average satisfies the
following constraint, which generalizes the 1D
constraint given by (\ref{eqn:roe_ave}):
\begin{equation}
\label{eqn:roe_ave_2d}
   \Phi^{\Tm} \equiv 
   \frac{1}{2}\sum_{i=1}^{3} \left( \vec{u}^{\, \Tm} \cdot \vec{n}_i \right) \, Q_i = 
   \oint_{\partial \Tm} \vec{f} (q^h) \cdot d\vec{s},
\end{equation}
where $q^h$ is the linear interpolant passing
through $(\vec{x}_i, Q_i)$ for $i=1,2,3$.
If the flux, $f(q)$, is at most a quadratic function of $q$,
then
\begin{equation}
   \vec{u}^{\, \Tm} = \frac{1}{3} \left( \vec{f}\,'(Q_1) + \vec{f}\,'(Q_2) + 
   \vec{f}\,'(Q_3) \right).
\end{equation}
The approximate Riemann solution gives rise to the
following set of fluctuations: 
\begin{equation}
\label{eqn:nscheme1}
  \text{\bf N-scheme:}  \quad
  \Phi^{\Tm}_i = \frac{1}{2} \left[ \vec{u}^{\, \Tm} \cdot \vec{n}_i \right]^{+}
  \left( Q_i - Q_{\star} \right), 
\end{equation}
where $Q_{\star}$ is the so-called {\it upwind parameter}.
In 1D the the upwind parameter relative to state $Q_i$ is
always either $Q_{i-1}$ if $u>0$ or $Q_{i+1}$ if $u<0$.
In multidimensions, $Q_{\star}$ is obtained by enforcing
the conservation constraint (\ref{eqn:cons_constraint}):
\begin{equation}
\label{eqn:nscheme2}
  Q_{\star} = \left( \sum_{i=1}^{3} [ \vec{u}^{\, \Tm} \cdot \vec{n}_i]^{-}
   Q_i
   \right) \Biggl/ \left( \sum_{j=1}^{3} [ \vec{u}^{\, \Tm} \cdot \vec{n}_j ]^{-}
   \right),
\end{equation}
where we have made use of the following two identities:
\begin{gather*}
   \sum_{i=1}^{3} \left( \vec{u}^{\, \Tm} \cdot \vec{n}_i \right) Q_i = 
   \sum_{i=1}^{3} \left[ \vec{u}^{\, \Tm} \cdot \vec{n}_i \right]^{+} Q_i 
   + \sum_{i=1}^{3} \left[ \vec{u}^{\, \Tm} \cdot \vec{n}_i \right]^{-} Q_i, \\
 \sum_{i=1}^{3} \left( \vec{u}^{\, \Tm} \cdot \vec{n}_i \right)  = 
   \sum_{i=1}^{3} \left[ \vec{u}^{\, \Tm} \cdot \vec{n}_i \right]^{+}  
   + \sum_{i=1}^{3} \left[ \vec{u}^{\, \Tm} \cdot \vec{n}_i \right]^{-} = 0.
\end{gather*}
In order to demonstrate that the N-scheme is monotonicity
preserving, we rewrite (\ref{eqn:nscheme1}) in the form (\ref{eqn:monotone_condition}) with
\begin{equation}
  c^{\Tm}_{ij} = \frac{\left[\vec{u}^{\, \Tm} \cdot \vec{n}_i \right]^{+}
  \left[\vec{u}^{\, \Tm} \cdot \vec{n}_j \right]^{-}}{\sum_{k}
  \left[\vec{u}^{\, \Tm} \cdot \vec{n}_k \right]^{-}} \ge 0.
\end{equation}
From this it can be shown that update (\ref{eqn:forwardeuler})
is monotone under the following CFL condition:
\begin{equation}
	\label{eqn:nscheme_cfl}
   \Delta t \le \min_i \left\{  \frac{2 | C_i |}{\sum_{\Tm: \, i \in \Tm}
   \left[ \vec{u}^{\, \Tm} \cdot \vec{n}_i \right]^{+}}  \right\}.
\end{equation}

\subsection{Linear preserving limiters}
The N-scheme described so far is both conservative
and monotonicity preserving; however, it is not yet
linear preserving. The problem with the previously
described N-scheme is that the weights (\ref{eqn:weights})
are not uniformly bounded independent of the mesh.
In order to modify the N-scheme to achieve uniformly
bounded $\beta_i$'s, Abgrall and Mezine \cite{article:AbMe04}
introduced a nonlinear limiting procedure. The limiting process
takes the original $\beta_i$ and replaces them with limited
versions, denoted $\tilde{\beta}_i$. The simpler of the
two approaches discussed in \cite{article:AbMe04}
yields the following formulas:
\begin{equation}
  \beta_i = \frac{\Phi^{\Tm}_i}{\Phi^{\Tm}} \quad \rightarrow \quad
  \tilde{\beta}_i = \frac{\left[ \beta_i \right]^{+}}{\sum_{j}
  \left[ \beta_j \right]^{+}},
\end{equation}
which guarantees that $0 \le \tilde{\beta}_i \le 1$.
The limited residuals are then given by
\begin{equation}
   \text{\bf Limited N-scheme:} \quad \Phi^{\Tm}_i = 
   \tilde{\beta}_i \, \Phi^{\Tm}.
\end{equation}
It is clear that this scheme is both linear preserving and
conservative. Furthermore, Abgrall and Mezine \cite{article:AbMe04} 
proved that the limited  N-scheme retains the monotonicity
properties of the original N-scheme with the same CFL
condition (\ref{eqn:nscheme_cfl}).

\subsection{Extension to systems}
\label{sec:rd_systems}
Following \cite{article:Decon93}, the N-scheme 
is extended to systems of conservation
laws by first defining the following averaged
flux Jacobians:
\begin{equation}
   \hat{J}^1 \equiv \frac{\partial f}{\partial q} \left( \hat{Z} \right),
   \quad \hat{J}^2 \equiv \frac{\partial g}{\partial q} \left( \hat{Z} \right),
   \quad \text{and} \quad
   \vec{J} \equiv \left( \hat{J}^1, \, \hat{J}^2 \right),
\end{equation}
where $Z$ is a parameterization of $Q$ and
\begin{equation}
   \hat{Z} = \frac{1}{3} \left( Z_1 + Z_2 + Z_3 \right).
\end{equation}
In order to achieve a conservative linearization
we must find a parameter vector, $Z$, such that
the following constraint is satisfied:
\begin{equation}
	\label{eqn:cons_nscheme_sys}
   \sum_{i=1}^3 \left( \vec{n}_i \cdot \vec{J} \right) \, \hat{Q}_i
   = \oint_{\partial \Tm} \vec{f} ( z^h ) \cdot d\vec{s},
\end{equation}
where $z^h$ is the linear interpolant that passes through
the points $(\vec{x}_i, Z_i)$ for $i=1,2,3$ and
\begin{equation}
    \hat{Q}_i \equiv \frac{\partial q}{\partial z}\left(\hat{Z}\right) \, Z_i.
\end{equation}
As was shown in \cite{article:Decon93} (see also
\cite{article:Ab01,article:Cs02}), constraint (\ref{eqn:cons_nscheme_sys}) will in general only be satisfied if
we are able to find a parameterization, $z = (z^1(q), 
z^2(q), \ldots, z^m(q))$, such that the flux, $\vec{f}(q(z))$,
depends at most quadratically on each $z^p$ for $p=1,2,\ldots,m$.
For the Euler equations from gas dynamics, as well as related
systems, such a parameterization is known \cite{article:Decon93}.

Assuming that a parameterization has been found,
we proceed by diagonalizing the flux Jacobian:
\begin{equation*}
   \vec{n}_i \cdot \vec{J} = R_i \Lambda_i R_i^{-1},
\end{equation*}
where $R_i$ is the matrix of right eigenvectors and $\Lambda_i$
is the diagonal matrix of eigenvalues.
Following the philosophy of Roe's approximate 
Riemann solver \cite{article:Roe81}, the systems N-scheme
is obtained by applying  the scalar N-scheme to each characteristic 
component. This results in the following method:
\begin{equation}
\label{eqn:sys_nscheme_1}
   \text{\bf Systems N-scheme:} \quad \Phi^{\Tm}_i
   = \frac{1}{2} R_i \Lambda^{+}_i R^{-1}_i \left( \hat{Q}_i - 
   {Q}_{\star} \right).
\end{equation}
The upwind parameter can be recovered by enforcing
local conservation:
\begin{equation}
  \label{eqn:sys_nscheme_2}
   Q_{\star} = \left\{ \sum_{i=1}^3 R_i \Lambda^{-}_i R^{-1}_i \right\}^{-1}
   \left\{  \sum_{j=1}^3 R_j \Lambda^{-}_j R^{-1}_j \, \hat{Q}_j \right\}.
\end{equation}
Note that solving for $Q_{\star}$ involves inverting
an $m \times m$ matrix.
Finally, we note that although it is based
on a generalization of the monotone scalar N-scheme, the
systems N-scheme is in general only approximately
non-oscillatory for nonlinear systems of conservation laws.
In practice, however, this scheme has been shown to work
quite well for steady-state shock computations for systems
such as the Euler equations from gas dynamics
\cite{article:Ab01}. 

The systems N-scheme  described so far 
is not linear preserving. In order that
the limiting procedure developed for scalar equations
can be re-used for systems,
Abgrall and Mezine \cite{article:AbMe04} proposed 
to project the distributed residuals into the eigenspace
of the Roe-Struijs-Deconinck-averaged flux Jacobian in some direction $\vec{n}$. In practice, the direction
$\vec{n}$ is chosen from physical considerations. For
example, in the case of the shallow water equations or the
Euler equation from gas dynamics, an approach that
gives good results in practice is to take $\vec{n}$ to be the
 local Roe-Struijs-Deconinck-averaged fluid velocity: $\vec{n} = \vec{u}$.
Once a limiting direction has been chosen, the limiting procedure can be summarized as follows:
\begin{align*}
   & \text{for } \, \,  p=1 \ldots m \\
   & \qquad \text{for } \, \, i=1,2,3:  \quad \text{set} \quad
     \Theta^p_i =  \ell^{p} \cdot \Phi^{\Tm}_i; \\
    & \qquad  \text{for } \, \, i=1,2,3:  \quad \text{set} \quad
    	\beta^p_i = \frac{\Theta^p_i}{\sum_{j=1}^{3} \Theta^p_j}; \\
   & \qquad  \text{for } \, \, i=1,2,3:  \quad \text{set} \quad
    	\tilde{\beta}^p_i = \frac{\left[ \beta^p_i \right]^{+}}{
	\sum_{j=1}^3 \left[ \beta^p_j \right]^{+}}; \\
   & \text{end} \\
   & \text{for } \, \, i=1,2,3:  \quad \text{set} \quad
   	 \Phi^{\Tm}_i = \sum_{p=1}^m  \tilde{\beta}_i^p \, \Theta^p_i
		\, r^p,
	\end{align*}
where $\ell^p$ and $r^p$ are the $p^{\text{th}}$ left and right 
eigenvectors of $\vec{n} \cdot \vec{J}$, respectively. 

\subsection{A correction for improved convergence}
\label{sec:Abfix}
As was pointed out by Abgrall \cite{article:Ab06}, the N-scheme
in conjunction with the limiting procedure outlined in
Section \ref{sec:rd_systems} has one major drawback:
the method does not in general converge to a steady-state
solution. The problem is not with the N-scheme itself, since
this method does converge to a steady-state, but instead
the problem lies in how the N-scheme interacts with the
limiting procedure. In the same paper, Abgrall  \cite{article:Ab06}
also provided a cure for this problem. 
He arrived at the following distributed residual:
\begin{equation}
   \Phi_i^{\Tm} = \underset{\text{limited N-scheme}}{\underbrace{B_i \, \Phi^{\Tm}}} +  \underset{\text{correction}}{\underbrace{
   \theta \, |\Tm|^{-1/2}
   \, K_i \, \Phi^{\Tm}}},
\end{equation}
where $K_i = \left( \vec{n}_i \cdot \vec{J} \right)/2$,
$|\Tm|$ is the area of triangle $\Tm$, and $\theta$ is a
grid and solution dependent parameter.
Notice that conservation is not affected by this correction term.
In order to produce a scheme that converges to a steady-state
solution, $\theta$ needs to be chosen so that the correction is
relatively small near shocks $(\theta = {\mathcal O}\left(|\Tm|^{1/2}\right))$ and relatively large in smooth regions $(\theta = 1)$.
Abgrall \cite{article:Ab06} proposed 
the following formula:
\begin{equation}
  \theta = \text{min} \left( 1, \frac{|\Tm|}{|\varphi^{\Tm}| + 
  10^{-10}} \right),
\end{equation}
where $\varphi^{\Tm}$ is the projection of $\Phi^{\Tm}$ onto
some important eigen-direction. In the case of the 
compressible Euler equations, $\varphi^{\Tm}$
should be taken to be the projection of $\Phi^{\Tm}$
onto the entropy wave.
In Section \ref{sec:num_examples}, in which we consider
several numerical examples, we will refer to the limited and
corrected N-scheme as the N$\backslash$LC-scheme.

\section{Multidimensional relaxation systems}
\label{sec:multid}
Having reviewed the relaxation scheme paradigm in
Section \ref{sec:review_1d} and residual distribution
schemes in Section \ref{sec:rd_schemes},
we now turn to develop a multidimensional relaxation system
framework. 
We again begin with the case of a scalar conservation law
and introduce the following relaxation system:
\begin{gather}
\label{eqn:relax_scalar_2d}
\begin{bmatrix}
  q \\ \mu^1 \\ \mu^2
\end{bmatrix}_{,t} + 
A^1
\begin{bmatrix}
  q \\ \mu^1 \\ \mu^2
\end{bmatrix}_{,x} +
A^2
\begin{bmatrix}
  q \\ \mu^1 \\ \mu^2
\end{bmatrix}_{,y} = \frac{1}{\varepsilon} 
\begin{bmatrix}
 0 \\
 f(q)  - \mu^1 \\
 g(q) - \mu^2
\end{bmatrix},
\end{gather}
where
\begin{align}
\label{eqn:relax_scalar_2d_A1}
A^1 &\equiv
\begin{bmatrix}
  0 \quad &  1 \quad & 0 \\
  - c^1 d^1 \quad & c^1 + d^1 \quad & 0 \\
  -c^2 d^1 - \frac{1}{4} \left( c^1 - d^1 \right) \left( c^2 - d^2 \right) \quad & \frac{1}{2} \left( c^2 + d^2 \right) \quad & \frac{1}{2} \left(
  	c^1 + d^1 \right)
\end{bmatrix}, \\
\label{eqn:relax_scalar_2d_A2}
A^2 &\equiv
\begin{bmatrix}
  0 \quad &  0 \quad & 1 \\
  -c^1 d^2 - \frac{1}{4} \left( c^1 - d^1 \right) \left( c^2 - d^2 \right) \quad & \frac{1}{2} \left( c^2 + d^2 \right) \quad & \frac{1}{2} \left(
  	c^1 + d^1 \right) \\
	- c^2 d^2 \quad & 0 \quad & c^2 + d^2 \quad
\end{bmatrix}.
\end{align}
In these expressions, we assume that $\vec{c}, \vec{d} \in \reals^2$.
Just as in the 1D case, we will separate the effects
of the hyperbolic left-hand side of this equation
from the relaxation source term on the right-hand
side by by viewing (\ref{eqn:relax_scalar_2d})
as being comprised of the following two sub-problems:
\begin{gather}
\label{eqn:relax_scalar_a}
\begin{bmatrix}
  q \\ \mu^1 \\ \mu^2
\end{bmatrix}_{,t} + 
A^1
\begin{bmatrix}
  q \\ \mu^1 \\ \mu^2
\end{bmatrix}_{,x} +
A^2
\begin{bmatrix}
  q \\ \mu^1 \\ \mu^2
\end{bmatrix}_{,y} = 0, \\
\label{eqn:relax_scalar_b}
\begin{bmatrix}
  \mu^1 \\ \mu^2
\end{bmatrix}_{,t}  = \frac{1}{\varepsilon} 
\begin{bmatrix}
 f(q)  - \mu^1 \\
 g(q) - \mu^2
\end{bmatrix}.
\end{gather}
In subsequent discussion we will make use of the following
matrix:
\begin{equation}
A(\vec{n}) = n^1 A^{1} + n^2 A^{2},
\end{equation}
where $\vec{n} \in \reals^2$ such that $\| \vec{n} \| \ne 0$.
The three eigenvalues of $A(\vec{n})$ are given by
\begin{align}
  \lambda^{1,3} &= \frac{1}{2} \left( \vec{c} \cdot \vec{n} +
  	\vec{d} \cdot \vec{n} \right) \pm \frac{1}{2}
	\sqrt{ \left(n^1 d^1 - n^1 c^1\right)^2 
	+ \left(n^2 d^2 - n^2 c^2 \right)^2}, \\
  \lambda^2 &= \frac{1}{2} \left( \vec{c} \cdot \vec{n} +
  	\vec{d} \cdot \vec{n} \right).
\end{align}
We note that these eigenvalues are strictly real.

Next we define a {\it relaxation procedure} that is
analogous to the 1D case; the main difference is that
in 2D a genuinely multidimensional Riemann
problem such as the one depicted in Figure \ref{fig:multid_riemann} must be solved. Instead of attempting to solve this exactly,
we solve it with the standard N-scheme. The full
procedure can then be summarized as follows:
\begin{enumerate}
  \item Choose values for the parameters $c^1$,
  $c^2$, $d^1$, and $d^2$.
  \item On an arbitrary triangle, ${\mathcal T}$, 
  	  approximately solve the multidimensional
	  Riemann problem associated with 
	  (\ref{eqn:relax_scalar_a}) by applying the standard 
	  N-scheme.
  \item Approximate the effect of equation
  	 (\ref{eqn:relax_scalar_b})
  	on the solution calculated in Step (2)
  	by directly setting  $\mu^1 = f(q)$
	and $\mu^2=g(q)$. In other words,	
 	instantaneously relax the solution from Step (2)
	to the $\varepsilon \rightarrow 0$ limit.
\end{enumerate}

\subsection{The N-scheme}
The first scheme that we will produce with the relaxation
procedure is the N-scheme
applied to the original scalar conservation law. 
We set
\begin{equation}
   \vec{c} = \vec{d} = \vec{u},
\end{equation}
where $\vec{u}$ is the Roe-Struijs-Deconinck-average
 that satisfies (\ref{eqn:roe_ave_2d}).
The above choice for $\vec{c}$ and $\vec{d}$
results in the following coefficient matrix for
the relaxation system:
\begin{equation}
  A(\vec{n}) = \begin{bmatrix}
  	0 \quad & n^1 \quad & n^2 \\
	-u^1 \left( \vec{u} \cdot \vec{n} \right) \quad &  n^1  u^1 +
	\vec{u} \cdot \vec{n} \quad &  n^2 u^1 \\
        -u^2 \left( \vec{u} \cdot \vec{n} \right) \quad & n^1 u^2 \quad & 
        n^2 u^2 + \vec{u} \cdot \vec{n}
  \end{bmatrix};
\end{equation}
this matrix has eigenvalues given by
\begin{equation}
  \lambda^{1,2,3} = \vec{u} \cdot \vec{n},
\end{equation}
and, as in the 1D case, has an incomplete set of
eigenvectors (i.e., the eigenvalues have
algebraic multiplicity 3, but geometric multiplicity of only 2).

Because $A(\vec{n})$ only has two linearly independent
eigenvectors it cannot be diagonalized;
and instead, we reduce it to Jordan canonical form via the following
similarity transformation:
\begin{equation}
    A(\vec{n}) = S \, M \, S^{-1},
\end{equation}
where 
\begin{align}
   M = \begin{bmatrix}
   	\vec{n} \cdot \vec{u} \quad & 1 \quad & 0 \\
	0 \quad & \vec{n} \cdot \vec{u} \quad & 0 \\
	0 \quad & 0 \quad & \vec{n} \cdot \vec{u}
   \end{bmatrix} \quad \text{and} \quad
   S = \begin{bmatrix}
   	1 \quad & -\frac{1}{\vec{n}\cdot\vec{u}} \quad & 0 \\
	u^1 \quad & -n^2 \quad & -n^2 \\
	u^2 \quad &  n^1 \quad &   n^1
   \end{bmatrix}.
\end{align}
In approximately solving the Riemann problem via
the N-scheme (Step (2) of the relaxation procedure),
we will need to make sense of the expression
$[ A(\vec{n}) ]^{+}$. Without a full set of eigenvectors,
we do this in the following way:
\begin{align}
   [ A(\vec{n}) ]^{+} \equiv
   S \, M^{+} \, S^{-1}, \quad
   \text{where } \, 
   M^{+} = \begin{bmatrix}
        \left[ \vec{n} \cdot \vec{u} \right]^{+} \quad & 1 \quad & 0 \\
	0 \quad & \left[ \vec{n} \cdot \vec{u} \right]^{+}  \quad & 0 \\
	0 \quad & 0 \quad & \left[ \vec{n} \cdot \vec{u} \right]^{+} 
   \end{bmatrix}.
\end{align}
This results in two possibilities:
\begin{enumerate}
   \item $\left[ \vec{n} \cdot \vec{u} \right]^{+} =  \left( 
   \vec{n} \cdot \vec{u} \right)
   \quad
   \Longrightarrow \quad \left[A(\vec{n})\right]^{+} = A(\vec{n})$;
   \item $\left[ \vec{n} \cdot \vec{u} \right]^{+}  = 0 \quad
   \Longrightarrow \quad \left[A(\vec{n})\right]^{+} 
   = \begin{bmatrix}
      - \left( \vec{n} \cdot \vec{u} \right)  \quad & n^1 \quad & n^2 \\
      - u^1 \left( \vec{n} \cdot \vec{u} \right) 
      \quad & u^1 n^1 \quad & u^1 n^2 \\
      - u^2 \left( \vec{n} \cdot \vec{u} \right) 
       \quad & u^2 n^1 \quad & u^2 n^2
   \end{bmatrix}$.
\end{enumerate}
The first of these two expressions is exactly the result
one should expect; the second expression, however, is
somewhat troubling. We should expect that
$\left[ A(\vec{n}_i) \right]^{+} \left( U_i - U_{\star} \right)
= 0$ if $\left( \vec{n}_i \cdot \vec{u}
\right) \le 0$, where
\begin{equation*}
  U_i \equiv \left( Q_i, \, \vec{f}(Q_i) \right) \quad \text{and} \quad
  U_{\star} \equiv  \left( Q_{\star}, \, \vec{\mu}_{\star} \right).
\end{equation*}
 Instead, we are currently stuck with the following result
 when $\left( \vec{n} \cdot \vec{u} \right) \le 0$:
\begin{equation*}
   \left[ A(\vec{n}_i) \right]^{+} \left( U_i - U_{\star} \right)
    = 
        \left\{ - \left( \vec{u}  \cdot \vec{n} \right) \left( Q_i - Q_{\star} \right)
         + \vec{n} \cdot \left( \vec{f}(Q_i) - \vec{\mu}^1_{\star} \right) 
         \right\}
         \begin{bmatrix}
         1 \\ u^1 \\ u^2
   \end{bmatrix}.
\end{equation*}
In order to clean up this result, we are forced to slightly 
modify the relaxation procedure for the N-scheme.
We will leave the sub-problem (\ref{eqn:relax_scalar_a}) alone,
but replace sub-problem (\ref{eqn:relax_scalar_b}) with the
the following system of ODEs:
\begin{equation}
\label{eqn:relax_scalar_bb}
\begin{bmatrix}
  \mu^1 \\ \mu^2
\end{bmatrix}_{,t}  = \frac{1}{\varepsilon} 
\begin{bmatrix}
 \frac{\partial f}{\partial q}\left( \hat{q} \right) \, q - \mu^1 \\
 \frac{\partial g}{\partial q}\left( \hat{q} \right) \, q - \mu^2
\end{bmatrix},
\end{equation}
where $\hat{q}$ is a piecewise constant function in space
that is constant on each triangle $\Tm$; the value
of this constant is  $\hat{Q}$, the multidimensional
Roe-Struijs-Deconinck average (\ref{eqn:roe_ave_2d}).
This results in the following modification of Step (3) in
the relaxation procedure:
\begin{enumerate}
   \setcounter{enumi}{2}
   \item On each triangle $\Tm$ approximate the effect of 
   (\ref{eqn:relax_scalar_bb})
   on the solution calculated in Step (2) by directly setting
   $\mu^1 = u^1 q$ and
   $\mu^2 = u^2 q$. In other
   words, instantaneously relax the solution from Step (2)
   to the $\varepsilon \rightarrow 0$ limit.
\end{enumerate}
Note that in general 
\[ \frac{\partial f}{\partial q}\left(\hat{q}\right) q \ne f(q) 
\quad \text{ and } \quad
\frac{\partial g}{\partial q}\left(\hat{q}\right) q \ne g(q);
\]
and therefore,
replacing (\ref{eqn:relax_scalar_b}) with (\ref{eqn:relax_scalar_bb})
will yield a different numerical scheme. As we will demonstrate
below, it is the scheme based on  (\ref{eqn:relax_scalar_bb})
that will reproduce the N-scheme.

The solution value at each node is now given by
\begin{equation}
  \label{eqn:nscheme_states}
   U_i \equiv \left( Q_i,  \, u^1  Q_i, \,  u^2 Q_i 
   \right).
\end{equation}
Additionally,  we enforce the condition:
\begin{equation}
	\label{eqn:nscheme_qstar}
   U_{\star} \equiv \left( Q_{\star}, \, u^1 Q_{\star}, \, u^2 Q_{\star}
   \right).
\end{equation}
With these modifications it is now true that
$\left[ A(\vec{n}_i) \right]^{+} \left( U_i - U_{\star} \right)
= 0$ if $\left( \vec{n}_i \cdot \vec{u}
\right) \le 0$.

In order to proceed with the relaxation procedure,
we solve a Riemann problem between
three states of the form (\ref{eqn:nscheme_states})
with $i=1,2,3$. Solving this Riemann problem via
the N-scheme tells us that the residuals distributed to each node
are given by
\begin{equation}
\label{eqn:relax_2d_res}
  \varphi^{\Tm}_i = \frac{1}{2} \left[ A\left( \vec{n}_i \right) \right]^{+}
  	 \left( U_i - U_{\star} \right) \equiv
	 \begin{cases}
	 	\frac{1}{2} A(\vec{n}_i) \left( U_i - U_{\star} \right) &
			\text{if} \quad \vec{u} \cdot \vec{n}_i > 0, \\
		0 & \text{otherwise},
	 \end{cases}
\end{equation}
where $U_{\star}$ has the form (\ref{eqn:nscheme_qstar}). 
In order to determine $U_{\star}$ in terms of the $U_i$
values, we add the three residuals given by (\ref{eqn:relax_2d_res})
and enforce that this sum yields the total residual in
triangle $\Tm$:
\begin{gather}
    \sum_{i=1}^3 \varphi^{\Tm}_i = \frac{1}{2} \sum_{i=1}^{3} 
    A(\vec{n}_i) \, U_i, \\
	\label{eqn:phitot_nscheme}
     \Longrightarrow   \quad
     \frac{1}{2} \sum_{i=1}^{3} \left[ A\left( \vec{n}_i \right) \right]^{+}
  	 \left( U_i - U_{\star} \right)  =
     \frac{1}{2} \sum_{i=1}^{3} \begin{bmatrix}
     	 \left( \vec{n}_i \cdot \vec{u} \right) Q_i \\
	 u^1 \, \left( \vec{n}_i \cdot \vec{u} \right) Q_i  \\
	 u^2 \, \left( \vec{n}_i \cdot \vec{u} \right) Q_i 
     \end{bmatrix}.
\end{gather}
If $\| \vec{u} \| > 0$, then we note the following result
on each triangle $\Tm$:
\begin{enumerate}
  \item $\exists \, k \in (1,2,3) \, \text{ such that }
  	\,  \vec{u} \cdot \vec{n}_k > 0$,
  \item $\exists \, k \in (1,2,3) \, \text{ such that }
  	\,  \vec{u} \cdot \vec{n}_k < 0$.
\end{enumerate}
This result implies that there are two possibilities whenever
$\| \vec{u} \| > 0$: the 1-target case -- $\exists$ exactly one
$k$ s.t. $\vec{u} \cdot \vec{n}_k > 0$, and the 2-target case --
$\exists$ exactly two $k$ s.t. $\vec{u} \cdot \vec{n}_k > 0$.
Without loss of generality, let us assume that 
$\vec{u} \cdot \vec{n}_1 > 0$ and $\vec{u} \cdot \vec{n}_3 < 0$, which yields one of the
two possibilities:
\begin{align}
\text{{\bf 1-target solution:}} & \quad 
\vec{u} \cdot \vec{n}_1 > 0, \quad \vec{u} \cdot \vec{n}_2 \le  0,
\quad \vec{u} \cdot \vec{n}_3 <  0, \\
\text{{\bf 2-target solution:}} & \quad
\vec{u} \cdot \vec{n}_1 > 0, \quad \vec{u} \cdot \vec{n}_2 > 0,
\quad \vec{u} \cdot \vec{n}_3 < 0.
\end{align}
We consider each of these two cases below.

\subsubsection{The 1-target solution}
The 1-target case is easy to analyze: the total residual
is completely distributed to the lone
node that is downwind of the flow, which we have taken without loss of generality to be node 1. If we let $\Phi$ denote the
component of the residual corresponding to $Q$, then
the 1-target case results in the following residual distribution:
\begin{equation}
   \Phi^{\Tm}_1 = \frac{1}{2} \sum_{i=1}^{3} \left( \vec{u} \cdot \vec{n}_i
   \right) Q_i,
   \quad \Phi^{\Tm}_2 = 0, \quad \Phi^{\Tm}_3 = 0,
\end{equation}
which is the same result that one would obtain
with the N-scheme on the original scalar conservation law.

\subsubsection{The 2-target solution}
The 2-target case involves distribution to two nodes, which
we have taken without loss of generality to be the 1 and 2 nodes.
From equation (\ref{eqn:phitot_nscheme})
we arrive at the following linear
system that must be solved in order to obtain
the upwind parameter $U_{\star}$:
\begin{equation}
- \left( A(\vec{n}_1) + A(\vec{n}_2)  \right)   U_{\star} =  A(\vec{n}_3) \, U_3.
\end{equation}
However, it is not difficult to show that 
\begin{equation}
  A(\vec{n}_1) + A(\vec{n}_2) + A(\vec{n}_3) = 0 
  \quad \Longrightarrow \quad  A(\vec{n}_3) =  -\left( A(\vec{n}_1) + A(\vec{n}_2)  \right).
\end{equation}
This implies that $U_{\star} \equiv U_3$. Therefore the 2-target
case results in the following residual distribution corresponding
to Q:
\begin{equation}
   \Phi^{\Tm}_1 = \frac{1}{2} \, \left( \vec{n}_1 \cdot \vec{u} \right) \left( Q_1 - Q_3 \right), \quad
    \Phi^{\Tm}_2 = \frac{1}{2} \, \left( \vec{n}_2 \cdot \vec{u} \right) \left( Q_1 - Q_3 \right),
   \quad \Phi^{\Tm}_3 = 0.
\end{equation}
This result is again identical to the original scalar conservation law.

\subsection{The RXN-scheme: genuinely multidimensional local
	Lax-Friedrichs}
One of the main difficulties with the N-scheme is that 
 computing the upwind parameter $Q_{\star}$ for
complicated systems of conservation laws can become
prohibitively expensive. Despite this fact, few alternatives
have been developed in the RD literature. One such
alternative was introduced by
Abgrall \cite{article:Ab06}, who 
considered a local Lax-Friedrichs-type method that was obtained, in analogy to the 1D case,
by taking the unstable ``centered''  residual and adding
the appropriate numerical viscosity.
In this work, we  construct 
a new method based on the idea of relaxation systems; this
scheme can
be viewed as a different multidimensional generalization
of the 1D LLF method.
For brevity we will call this method the RXN-scheme,
which stands for ``relaxation N-scheme\footnote{the words
``N-scheme'' appear here because we make use of the
N-scheme to solve the homogeneous part of the relaxation
system.}.'' In analogy with
the 1D LLF method as derived in Section \ref{sec:LLF1D},
we make the following choice for the 
parameters $\vec{c}$ and $\vec{d}$ in 
(\ref{eqn:relax_scalar_2d})--(\ref{eqn:relax_scalar_2d_A2}):
\begin{equation}
  \vec{d} = - \vec{c} =  \bigl( s^{\Tm}, \, s^{\Tm} \bigr).
\end{equation} 
Note that each triangle can have a different value of
$s^{\Tm}$; this is why we call it a `local' Lax-Friedrichs.
This choice results in a coefficient matrix of the form
\begin{equation}
	\label{eqn:coeff_RXN}
     A(\vec{n}) = \begin{bmatrix}
     0 \quad & n^1 \quad & n^2 \\
     n^1 \, \left( s^{\Tm} \right)^2 \quad & 0 \quad & 0 \\
     n^2 \, \left( s^{\Tm} \right)^2 \quad & 0 \quad & 0 
     \end{bmatrix},
\end{equation}
which has an eigenvector decomposition given by
\begin{equation}
\begin{split}
	&A(\vec{n}) = S M S^{-1}  \\
	&=\begin{bmatrix}
	-1 \quad & 0 \quad & 1 \\
	\frac{s^{\Tm} n^1}{\| \vec{n} \|} \quad & -\frac{n^2}{\| \vec{n} \|}
	 \quad & \frac{s^{\Tm} n^1}{\| \vec{n} \|} \\
	\frac{s^{\Tm} n^2}{\| \vec{n} \|} \quad & \frac{n^1}{\| \vec{n} \|}
	 \quad & \frac{s^{\Tm} n^2}{\| \vec{n} \|}
	\end{bmatrix}
	\begin{bmatrix}
	-\| \vec{n} \|  s^{\Tm} & & \\ & 0 &  \\  & & \| \vec{n} \|  s^{\Tm}
	\end{bmatrix}
	\begin{bmatrix}
		-\frac{1}{2} \quad & \frac{n^1}{2 \| \vec{n} \| s^{\Tm}} \quad & 
		\frac{n^2}{2 \| \vec{n} \| s^{\Tm}} \\
		 0 \quad & -\frac{n^2}{\| \vec{n} \|} \quad & \frac{n^1}{\| \vec{n} \|} \\
		 \frac{1}{2} \quad & \frac{n^1}{2 \| \vec{n} \| s^{\Tm}} \quad & 
		\frac{n^2}{2 \| \vec{n} \| s^{\Tm}}
	\end{bmatrix}.
\end{split}
\end{equation}
The Chapman-Enskog expansion for this relaxation system
can be written to first order as
\begin{equation}
   q_{,t} + f(q)_{,x} + g(q)_{,y} \approx \varepsilon \, \nabla \cdot \left(  
   	\begin{bmatrix}
	  \left(s^{\Tm}\right)^2 - \left( f'(q) \right)^2 \quad &  -f'(q) \, g'(q) \\
	   -f'(q) \, g'(q) \quad & \left(s^{\Tm}\right)^2 - \left( g'(q) \right)^2
	\end{bmatrix}
	\nabla q \right).
\end{equation}
The eigenvalues of the diffusion matrix in the above expression
are 
\begin{equation}
	\lambda^1 = \left(s^{\Tm}\right)^2, \quad \lambda^2 = 
	\left(s^{\Tm}\right)^2 - \left( f'(q)^2 + g'(q)^2 \right ),
\end{equation}
which results in the following restriction on the choice of the
lone parameter $s^{\Tm}$:
\begin{equation}
\label{eqn:subchar_rxn}
  s^{\Tm} \ge \| \vec{f}\,'(q) \|,
\end{equation}
for all $q \in \Tm$.

Applying the N-scheme to the relaxation system
with coefficient matrix (\ref{eqn:coeff_RXN}), yields the following
residual
\begin{equation}
  \varphi_i  =  
  \frac{1}{2} S_i M_i^{+} S_i^{-1} \left(  U_i - U_{\star}
  \right), 
\end{equation}
where 
$U_{\star} = \left( Q_{\star}, \, \mu_{\star}^1, \, \mu_{\star}^2 \right)$.
Simplifying this expression gives
\begin{equation}
\varphi_i = \frac{1}{4} \biggl\{  s^{\Tm} \| \vec{n}_i \| \left( Q_i - Q_{\star} \right)
	 + \vec{n}_i \cdot \left( \vec{f}(Q_i) - \vec{\mu}_{\star} \right) 
	 \biggr\}
	\begin{bmatrix}
	  1 \\ s^{\Tm} \frac{n^1_i}{\| \vec{n}_i \|} \\ s^{\Tm} \frac{n^2_i}{\| \vec{n}_i \|}
	\end{bmatrix}.
\end{equation}
In order to calculate $U_{\star}$, we must enforce
conservation:
\begin{equation}
\sum_{i=1}^{3} \varphi_i = \varphi^{\Tm} = 
\frac{1}{2} \sum_{i=1}^{3} 
\begin{bmatrix}
	\vec{n}_i \cdot \vec{f}(Q_i) \\
	n^1_i \, \left( s^{\Tm} \right)^2 \, Q_i \\ n^2_i \,  \left(
	s^{\Tm} \right)^2 \, Q_i 
\end{bmatrix},
\end{equation}
which results in the following linear system 
for the upwind parameters $(Q_{\star}, \vec{\mu}_{\star})$:
\begin{equation}
\label{eqn:rxn_star_system}
\sum_{i=1}^{3}
  \begin{bmatrix}
  	s^{\Tm} \| \vec{n}_i \| \quad & 0 \quad & 0 \\
	0 \quad &  \frac{n^1_i \, n^1_i}{\| \vec{n}_i \|} \quad &
	 \frac{n^1_i \, n^2_i}{\| \vec{n}_i \|} \\
	0 \quad &  \frac{n^1_i \, n^2_i}{\| \vec{n}_i \|} \quad &
	 \frac{n^2_i \, n^2_i}{\| \vec{n}_i \|}
  \end{bmatrix} 
  \begin{bmatrix}
  	Q_{\star} \\ \mu^1_{\star} \\ \mu^2_{\star}
  \end{bmatrix} = 
  \sum_{j=1}^{3}
  \begin{bmatrix}
  	s^{\Tm} \| \vec{n}_j \| Q_j - \vec{n}_j \cdot \vec{f}(Q_j) \\
	 n^1_j \left( \frac{\vec{n}_j}{\| \vec{n}_j \|}  \cdot \vec{f}(Q_j) - s^{\Tm} Q_j \right) \\
         n^2_j \left( \frac{\vec{n}_j}{\| \vec{n}_j \|}  \cdot \vec{f}(Q_j) - s^{\Tm} Q_j \right) 
  \end{bmatrix}.
\end{equation}
The solution to this linear system can be written as
\begin{align}
	\label{eqn:qstar_rxn}
   Q_{\star} &= \frac{\sum_{j=1}^3 \left( s^{\Tm} \| \vec{n}_j \| Q_j - \vec{n}_j
   	\cdot \vec{f}(Q_j) \right) }{\sum_{i=1}^3 s^{\Tm} \| \vec{n}_i \|}, \\
	\label{eqn:mustar1_rxn}
   \mu^1_{\star} &=  \frac{1}{N} \sum_{i=1}^3 \sum_{j=1}^3 
   	\frac{n^2_i}{\| \vec{n}_i \|}
	\left(  \frac{\vec{n}_j}{\| \vec{n}_j \|} \cdot \vec{f}(Q_j)
	- s^{\Tm} Q_j \right)
	 \left(  n^2_i n^1_j - n^2_j n^1_i  \right), \\
	 \label{eqn:mustar2_rxn}
	\mu^2_{\star} &=  \frac{1}{N} \sum_{i=1}^3 \sum_{j=1}^3 
   	\frac{n^1_i}{\| \vec{n}_i \|} \left(  \frac{\vec{n}_j}{\| \vec{n}_j \|} \cdot \vec{f}(Q_j)
	- s^{\Tm} Q_j \right) \left(  n^1_i n^2_j - n^1_j n^2_i  \right),
\end{align}
where
\begin{equation}
  \label{eqn:mustar3_rxn}
   N = \sum_{i=1}^3 \sum_{j=1}^3 \left\{   \frac{n^1_i \,
   n^1_i \, n^2_j \, n^2_j - n^1_i \, n^1_j \, n^2_i \, n^2_j}{\| \vec{n}_i \| 
   \, \| \vec{n}_j \|}   \right\}.
\end{equation}

Let us now take a moment to reflect on what just happened.
Although the original coefficient matrix, (\ref{eqn:coeff_RXN}), for this
method was comically simple, after applying the N-scheme to this
system on an arbitrary triangle, the resulting upwind parameters
are somewhat complicated. On the other hand, we see
from equation (\ref{eqn:rxn_star_system}) that the parameter
$Q_{\star}$ is completely decoupled from $\vec{\mu}_{\star}$.
We make use of this last fact to construct an alternative 
scheme in the following way: instead of computing the
components of $\vec{\mu}_{\star}$ from (\ref{eqn:mustar1_rxn})--(\ref{eqn:mustar3_rxn}),
we enforce 
\begin{equation}
     \vec{\mu}_{\star} \equiv \vec{f}(Q_{\star})
\end{equation}
by again invoking the $\varepsilon \rightarrow 0$ relaxation limit.
Although this direct enforcement clearly gives a different
scheme than if we had used (\ref{eqn:mustar1_rxn})--(\ref{eqn:mustar3_rxn}), what we achieve with this approach
is a very simple method that we refer to as the RXN-scheme
(relaxation N-scheme). In terms of the residual distributed to
node $i$ in the $Q$-variable, we now obtain the following expression:
\begin{equation}
\label{eqn:RXN_scheme}
\text{\bf RXN-scheme:} \quad
   \Phi^{\Tm}_i  =  \frac{1}{4} \, s^{\Tm} \, \| \vec{n}_i \| \left( Q_i - Q_{\star} \right)
   	+  \frac{1}{4} \, \vec{n}_i \cdot \left( \vec{f}(Q_i) - \vec{f}(Q_{\star})  \right),
\end{equation}
where $Q_{\star}$ is given by (\ref{eqn:qstar_rxn}).
Note that this method is automatically conservative since 
$Q_{\star}$  still satisfies the first equation in linear
system (\ref{eqn:rxn_star_system}).

\begin{theorem}
\emph{(Monotonicity)}
If there exists a $\tilde{Q}$ such that
\begin{equation}
\label{eqn:thm1}
   \sum_{i=1}^{3} \vec{n}_i \cdot \vec{f}(Q_i) = 
    \sum_{i=1}^{3} \vec{n}_i \cdot \vec{f}'(\tilde{Q}) \, Q_i, 
\end{equation}
and
\begin{equation}
\label{eqn:thm2}
	s^{\Tm} \ge \max \left\{ \| f'(\tilde{Q}) \|, \,
	  \| f'(\bar{Q}_1) \|, \, \| f'(\bar{Q}_2) \|, \,
	  \| f'(\bar{Q}_3) \| \right\},
\end{equation}
where for each $j$
\begin{equation}
	\label{eqn:thm3}
    \vec{n}_j \cdot \left( \vec{f}(Q_j) - \vec{f}(Q_{\star}) \right)
    = \vec{n}_j \cdot \vec{f}\,' \left( \bar{Q}_j \right)
    \left( Q_j - Q_{\star} \right)
\end{equation}
from the Rankine-Hugoniot conditions, then
the RXN-scheme as defined by (\ref{eqn:RXN_scheme})
and (\ref{eqn:qstar_rxn}) satisfies
the following condition:
\begin{equation}
\label{eqn:monotone}
   \Phi^{\Tm}_i = \sum_{j=1}^{3} c^{\Tm}_{ij} \left( Q_i - Q_j \right),
\end{equation}
where $c^{\Tm}_{ij} \ge 0$ for all $i,j=1,2,3$. This condition
along with the following CFL constraint on the 
time step:
\begin{equation}
  \label{eqn:RXN_timestep1}
   \Delta t \, \le \,
	 \min_{i} \left\{\frac{2 |C_i|}{\sum_{\Tm: i \in \Tm}
   	\| \vec{n}_i \| \, s^{\Tm}}  \right\},
\end{equation}
is enough to guarantee that the RXN-scheme with
forward Euler time-stepping (\ref{eqn:forwardeuler}) is monotonicity
preserving.
\end{theorem}

\begin{proof}
(1) Using the Rankine-Hugoniot conditions (\ref{eqn:thm3}),
we rewrite the RXN-scheme as
\begin{equation}
\label{eqn:proof_part1}
   \Phi_i = P^{\Tm}_i \left( Q_i - Q_{\star} \right), \quad \text{where}
   \quad P^{\Tm}_i \equiv \frac{1}{4} \left( s^{\Tm} \| \vec{n}_i \| + \vec{n}_i \cdot
   \vec{f}\,' \left( \bar{Q}_i \right) \right). 
\end{equation}
Similarly, we rewrite (\ref{eqn:qstar_rxn}) as follows
\begin{equation}
\label{eqn:proof_part2}
    Q_{\star} = \frac{\sum_{j=1}^3 N^{\Tm}_j \, Q_j}{\sum_{j=1}^3 N^{\Tm}_j},
     \quad \text{where} \quad
    N^{\Tm}_j \equiv s^{\Tm} \| \vec{n}_j \| - \vec{n}_j \cdot \vec{f}\,'(\tilde{Q}).
    \end{equation}
    Note that the above expression was obtained by making
    use of (\ref{eqn:thm1}) and 
    the identity: $\sum_{k=1}^3 \vec{n}_k \cdot
    	\vec{f}\,'(\hat{Q}) = 0$.
Combining expressions (\ref{eqn:proof_part1}) and (\ref{eqn:proof_part2}) yields (\ref{eqn:monotone})
with
\begin{equation}
   c^{\Tm}_{ij} \equiv \frac{P^{\Tm}_i \, N^{\Tm}_j}{\sum_{k=1}^3 N^{\Tm}_k}.
\end{equation}
We note that $c^{\Tm}_{ij} \ge 0$ $\forall i,j \in (1,2,3)$, 
because (\ref{eqn:thm2}) implies that
$P^{\Tm}_i \ge 0$ $\forall i \in (1,2,3)$ and $N^{\Tm}_j \ge 0$ $\forall j \in (1,2,3)$.

(2) We now insert expression (\ref{eqn:monotone}) into
(\ref{eqn:forwardeuler}) and simplify:
\begin{align*}
   & Q^{n+1}_i  =  Q^n_i - \frac{\Delta t}{| C_i |}
   \sum_{\Tm: \, i \in \Tm} \, \sum_{j \in \Tm} \,
   c_{ij}^{\Tm} \left( Q^n_i - Q^n_j \right), \\
   \Longrightarrow \quad & Q^{n+1}_i  =  \left\{ 1 - 
   \frac{\Delta t}{| C_i |} \, \sum_{\Tm: \, i \in \Tm} \, 
   \sum_{j \in \Tm} c^{\Tm}_{ij}  \, 
    \right\} Q^n_i + \frac{\Delta t}{| C_i |} \, \sum_{\Tm: \, i \in \Tm} \, 
   \sum_{j \in \Tm}  c_{ij}^{\Tm} \, Q^n_j.
\end{align*}
Monotonicity is achieved if $Q^{n+1}_i$ is a convex
average of all of the surrounding $Q^n_j$. Since
each $c^{\Tm}_{ij} \ge 0$, we obtain a convex average
provided that 
\begin{gather*}
  \frac{\Delta t}{| C_i |} \, \left( \sum_{\Tm: \, i \in \Tm} \, 
   \sum_{j \in \Tm} c^{\Tm}_{ij} \right) \, \le 1
   \quad \Longrightarrow \quad
   \Delta t_i \le \frac{| C_i |}{\left( \sum_{\Tm} \, 
   \sum_{j} c^{\Tm}_{ij} \right)}  
   = 
    \frac{| C_i |}{\left( \sum_{\Tm} \, 
   \sum_{j} \frac{P^{\Tm}_i N^{\Tm}_j}{\sum_{k} N^{\Tm}_k} 
   \right)}, \\
  \Longrightarrow \quad \Delta t \, \le \, \min_{i} \left\{  \frac{|C_i|}{\sum_{\Tm: \, i \in \Tm}
   	P_i^{\Tm}}  \right\}.
\end{gather*}
The  time restriction is clearly satisfied if we take
(\ref{eqn:RXN_timestep1}).
\qed
\end{proof}

In practice the time step presented in the above
theorem is overly restrictive. In the numerical
simulations presented in Section \ref{sec:num_examples},
we instead use the same time-step as used with the N-scheme:
$85\%$ of the maximum CFL number given by
expression (\ref{eqn:nscheme_cfl}).

\subsection{Systems N-scheme}
The systems generalization of coefficient matrix
(\ref{eqn:coeff_RXN}) for a system of $m$ conserved
variables is the following $3m \times 3m$ matrix:
\begin{equation}
   A(\vec{n}) = 
   \begin{bmatrix}
   	0 \id \quad & n^1 \, \id  \quad & n^2  \, \id \\
	- \left(\vec{n} \cdot \vec{J}\right) \hat{J}^1 \quad
	&  n^1 \hat{J}^1 + \vec{n} \cdot \vec{J} \quad &
	n^2 \hat{J}^1 \\
	- \left(\vec{n} \cdot \vec{J}\right) \hat{J}^2 \quad
	& n^1 \hat{J}^2 &  n^2 \hat{J}^2 + \vec{n} \cdot \vec{J} \quad
	 \\
   \end{bmatrix},
\end{equation}
where $\id$ is again the $m \times m$ identity matrix,
$0 \id$ is the $m \times m$ matrix with zeros in
every entry, and
$\vec{J} = \left( \hat{J}^1, \, \hat{J}^2 \right)^t$ is
the flux Jacobian matrix evaluated at the Roe-Struijs-Deconinck
average \cite{article:Decon93}.
The systems generalization of 
(\ref{eqn:nscheme_states})--(\ref{eqn:nscheme_qstar})
are the following $3 m \times 1$ vectors:
\begin{equation}
  U_i \equiv \left( Q_i, \, \hat{J}^1 Q_i, \, \hat{J}^2 Q_i \right) 
  \quad \text{and} \quad
  U_{\star} \equiv \left( Q_{\star}, \, \hat{J}^1 Q_{\star}, \, \hat{J}^2 Q_{\star} \right).
\end{equation}
In order to calculate the appropriate residuals in the
relaxation procedure, we need to again understand how
to create the matrices $\left[A (\vec{n}) \right]^{+}$ and 
$\left[A (\vec{n}) \right]^{-}$.
As in the scalar case this is complicated by the fact that
$A(\vec{n})$ does not have a full set of eigenvectors. In particular,
the Jordan canonical form of this matrix can be written as
\begin{equation}
   A(\vec{n}) = S 
   \begin{bmatrix}
   D^1    \\
   \quad &    \ddots \\
   \quad & \quad &
   D^m
   \end{bmatrix}
   S^{-1}, \quad
   \text{where} \quad 
    D^p = \begin{bmatrix}
   \lambda^p \quad & 1 \quad & 0 \\
   0 \quad & \lambda^p \quad & 0 \\
   0 \quad & 0 \quad & \lambda^p
   \end{bmatrix}.
\end{equation}
Here $\lambda^p$ is the $p^{\text{th}}$ eigenvalue
of the $m \times m$ matrix $\vec{n} \cdot \vec{J}$.
We omit the complicated expression for the matrix $S$.
In order to obtain an expression for $\left[A (\vec{n}) \right]^{+}$, one
has to replace each $\lambda^p$ in the above expression
with $\left[ \lambda^p \right]^{+}$. Carrying this out
results in the following matrix:
\begin{equation}
   \left[A (\vec{n}) \right]^{+} = 
   \begin{bmatrix}
   	- \hat{J}^{-} \quad &
	  n^1 \, \id  \quad & n^2  \, \id \\
	- \hat{J}^{+} \hat{J}^1
	- \hat{J}^1 \hat{J}^{-}  \quad
	&  n^1 \hat{J}^1 + \hat{J}^{+} \quad &
	n^2 \hat{J}^1 \\
	- \hat{J}^{+} \hat{J}^2
	- \hat{J}^2 \hat{J}^{-} \quad
	& n^1 \hat{J}^2 &  n^2 \hat{J}^2 + \hat{J}^{+} 
   \end{bmatrix},
\end{equation}
where $\hat{J}^{\pm} = \bigl( \vec{n} \cdot \vec{J}
\bigr)^{\pm} = R \Lambda^{\pm} R^{-1}$, $\Lambda$ is
the diagonal matrix of eigenvalues of $\vec{n} \cdot \vec{J}$,
and $R$ is the corresponding matrix right-eigenvectors.
An analogous formula for $\left[A (\vec{n}) \right]^{-}$ can also be
readily constructed.
From the above expression we find that
\begin{equation}
  \left[A (\vec{n}) \right]^{\pm} \, U_i 
  = \left[ \left( \vec{n} \cdot \hat{J} \right)^{\pm} Q_i, \, \,
   \hat{J}^1 \left( \vec{n} \cdot \hat{J} \right)^{\pm} Q_i, \, \,
   \hat{J}^2 \left( \vec{n} \cdot \hat{J} \right)^{\pm} Q_i \right]^{t}.
\end{equation}

Having established expressions for $\left[A (\vec{n}) \right]^{\pm}$, we now
 proceed by applying the N-scheme to the relaxation system:
\begin{gather}
  \label{eqn:linsys_nrelax}
   \varphi_i = \frac{1}{2} \left[ A(\vec{n}_i) \right]^{+} \left( U_i - U_{\star} \right) 
   \,
   \Longrightarrow \,
      \sum_{i=1}^3 \left[A (\vec{n}_i) \right]^{-} 
      \, U_{i} =  \left( - \sum_{i=1}^3 \left[A (\vec{n}) \right]^{+} \right) \, U_{\star},
\end{gather}
where
\begin{equation}
  - \sum_{i=1}^3 \left[A (\vec{n}_i) \right]^{+} = 
  \begin{bmatrix}
     \sum_i \hat{J}_i^{-} \quad &
     0 \quad & 0 \\
     - \sum_i \left[ \hat{J}_i^{-}  \hat{J}^1
     +  \hat{J}^1 \hat{J}_i^{-}  \right]
      \quad &
     \sum_i \hat{J}_i^{-} \quad & 0 \\
      - \sum_i \left[ \hat{J}_i^{-}  \hat{J}^2
     +  \hat{J}^2 \hat{J}_i^{-}  \right]
 	\quad & 0 \quad &
     \sum_i \hat{J}_i^{-}
  \end{bmatrix}
\end{equation}
and $\hat{J}_i = \left( \vec{n}_i \cdot \vec{J} \right)$.
The unique solution to the linear system in (\ref{eqn:linsys_nrelax}) is
$U_{\star} = \left( Q_{\star}, \, \hat{J}^1 Q_{\star}, \, \hat{J}^2 Q_{\star}
\right)$ with
\begin{equation}
   Q_{\star} = \left\{ \sum^3_{i=1} \left( \vec{n}_i \cdot 
   \vec{J} \right)^{-} \right\}^{-1} \left\{ \sum^3_{i=1} \left( \vec{n}_i \cdot 
   \vec{J} \right)^{-} Q_i  \right\},
\end{equation}
and the component of the residual $\varphi_i$
associated with $Q$ can be written as
\begin{equation}
   \Phi_i^{\Tm} = \frac{1}{2} \left( \vec{n}_i \cdot 
   \vec{J} \right)^{+} \left( Q_i - Q_{\star} \right). 
\end{equation}
This result shows that this relaxation scheme identically reproduces 
the systems N-scheme (\ref{eqn:sys_nscheme_1})--(\ref{eqn:sys_nscheme_2}).

\subsection{Systems RXN-scheme}
Just as the LLF method in the one-dimensional case, the RXN-scheme
extends to systems of conservation laws in a simple manner. All that
we have to do is apply the scalar version of the scheme to each
component of the vector conserved variable. The coefficient
matrix in the relaxation procedure can be written as
\begin{equation}
A(\vec{n}) = 
\begin{bmatrix}
   0 \id \quad & n^1 \, \id  \quad & n^2  \, \id \\
   n^1 \, s^2 \, \id \quad & 0  \, \id \quad & 0  \, \id \\
   n^2 \, s^2 \, \id \quad & 0  \, \id \quad & 0  \, \id
\end{bmatrix},
\end{equation}
where $\id$ is the $m \times m$ identity matrix.
In order to satisfy the sub-characteristic condition
we require that
\begin{equation}
   s \ge \max_{p=1, \ldots, m} \sqrt{ \left[ \lambda^{p,x}(q) \right]^2
   + \left[ \lambda^{p,y}(q) \right]^2},
\end{equation}
over all $q \in \Tm$. In the above expression
$\lambda^{p,x}$ and $\lambda^{p,y}$ are the $p^{\text{th}}$ eigenvalue of $\partial f/\partial q$
and $\partial g/\partial q$, respectively.

We note that the systems RXN-scheme, and in particular, the
version of this scheme with limiters (Section \ref{sec:rd_systems})
and convergence corrections (Section \ref{sec:Abfix}),
provides an alternative to the systems N-scheme that does not require
the inversion of an $m \times m$ matrix in each element at each
time level, nor does it require any special entropy fixes or 
special treatment near stagnation points. Since this method is
also simpler than the N-scheme, it should also yield some gains
in computational efficiency. The systems N-scheme and
RXN-scheme are compared in detail in Section 
\ref{sec:num_examples}.

\subsection{RXN-scheme in $d$-dimensions}
The above procedure for obtaining the 2D RXN-scheme can
be generalized to any space dimension. In the $d$-dimensional 
case we arrive at the following scheme:
\begin{equation}
\text{\bf RXN}_d \text{\bf \, -scheme:} \quad
   \Phi^{\Tm}_i  =  \frac{s^{\Tm} \, \| \vec{n}_i \| \left( Q_i - Q_{\star} \right)
   	+  \vec{n}_i \cdot \left( \vec{f}(Q_i) - \vec{f}(Q_{\star})  \right)}{2d},
\end{equation}
where 
\begin{equation}
Q_{\star} = \frac{\sum_{j=1}^{d+1} \left( s^{\Tm} \| \vec{n}_j \| Q_j - \vec{n}_j
   	\cdot \vec{f}(Q_j) \right) }{\sum_{i=1}^{d+1} s^{\Tm} \, \| \vec{n}_i \|}.
\end{equation}
In particular, we note that for $d=1$, this scheme exactly
reduces to the 1D local Lax-Friedrichs method 
\cite{article:Ru61}. We also mote that the $1/d$ geometric factor
comes from the $d$-dimensional N-scheme; see 
for example equation (7) in \cite{article:Cs02}.

\section{Numerical examples}
\label{sec:num_examples}
In this section we compare the N-scheme and
the newly proposed RXN scheme on several
numerical examples. We will refer to the versions of
the N-scheme and RXN-scheme that have
been limited according to the procedure outlined in Section
\ref{sec:rd_systems} and corrected according to the procedure outlined in 
Section \ref{sec:Abfix} as the N$\backslash$LC-scheme and
RXN$\backslash$LC-scheme, respectively.

\subsection{Steady-state advection}
First, we consider the advection equation on 
$[0,1] \times [0,1]$:
\begin{equation}
    q_{,t} + \vec{u} \cdot \nabla q = 0,
\end{equation}
with non-divergent velocity and boundary conditions
given by
\begin{gather*}
	\vec{u}(x,y) = \left( -\pi y, \, \pi x \right), \\
	q(1,y) = 0, \quad 
	q(x,0) = \begin{cases}
		\sin \left( \pi \left( \frac{0.7 - x }{0.6} \right) \right)
		& \quad \text{if \, } 0.1 < x < 0.7, \\
		0 & \quad \text{otherwise}.
	\end{cases}
\end{gather*}
This same problem was considered in \cite{article:Ab06}.

For a non-divergent velocity field, an elegant way to
solve the advection equation using the N-scheme
is through the introduction of a streamfunction:
\begin{align*}
	\psi(x,y) = -\frac{\pi}{2} \left( x^2 + y^2 \right),
\end{align*}
such that $\vec{u} = \left( \partial \psi/\partial y, \, -\partial \psi/\partial x
\right)$.
The N-scheme can then be written as
\begin{equation*}
    \Phi^T_i  =  k_i^{+} \left( Q_i  - Q_{\star} \right),
\end{equation*}
where
\begin{align*}
    k_1  &=  \frac{1}{2} \left( \psi(x_2, y_2) - \psi(x_3, y_3) \right), \\
    k_2  &=  \frac{1}{2} \left( \psi(x_3, y_3) - \psi(x_1, y_1) \right), \\
    k_3  &=  \frac{1}{2} \left( \psi(x_1, y_1) - \psi(x_2, y_2) \right).
\end{align*}
The advantage of this formulation is that we
achieve numerical conservation,  even though the equations are solved in advective form.

For the RXN scheme, we use residual (\ref{eqn:RXN_scheme}) where
 the flux functions are given by
\[
	\vec{f}(Q_i) = \vec{u}_i  Q_i, \quad  \vec{f}(Q_{\star}) 
	= \vec{u}_{\star} Q_{\star},
	\quad \vec{u}_{\star} = 
	\frac{\| \vec{n}_1 \| \vec{u}_1 + \| \vec{n}_2 \| \vec{u}_2  
	+ \| \vec{n}_3 \| \vec{u}_3}{\| \vec{n}_1 \| + \| \vec{n}_2 \| + \| \vec{n}_3 \|},
\]
and $Q_{\star}$ is given by (\ref{eqn:rxn_star_system}) as usual.

\begin{figure}[!t]
\begin{center}
\mbox{(a)  \includegraphics[height=62mm]{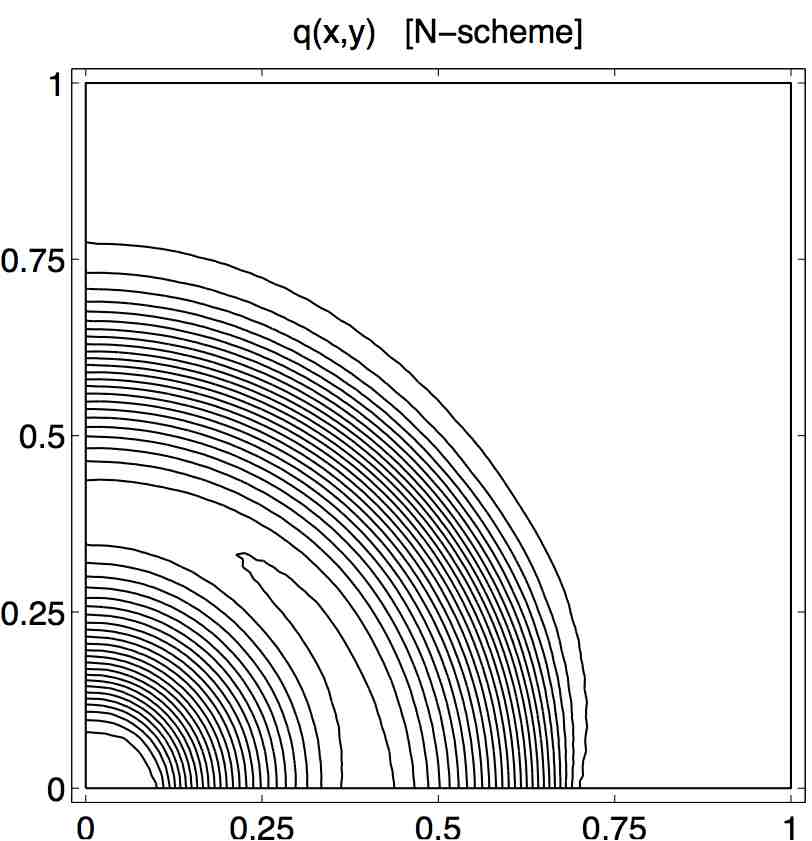} \,
(b)  \includegraphics[height=62mm]{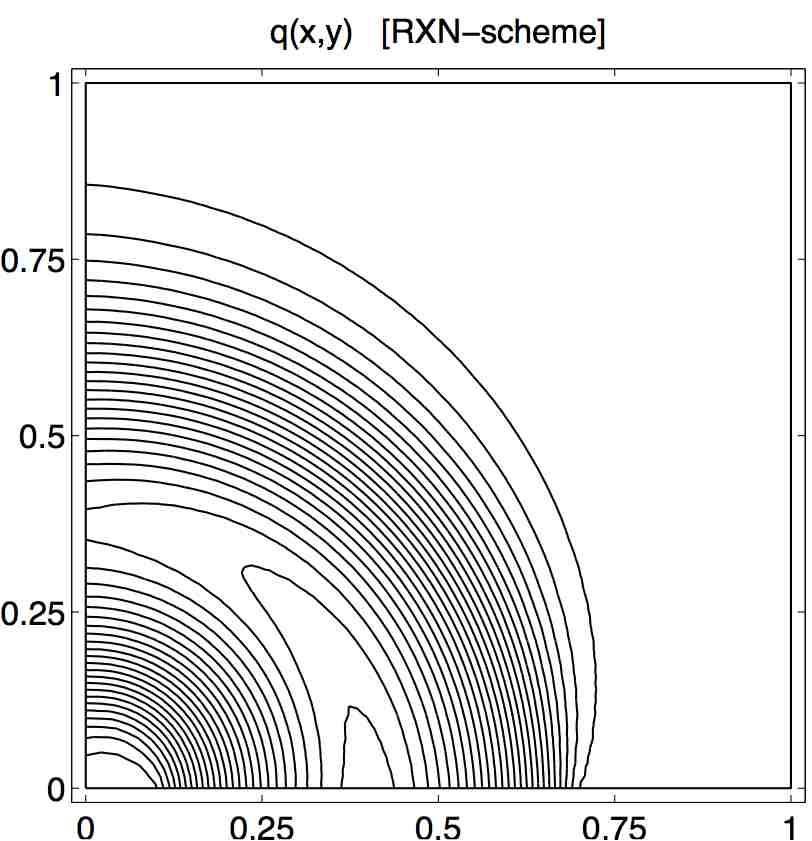}}

\mbox{(c)  \includegraphics[height=62mm]{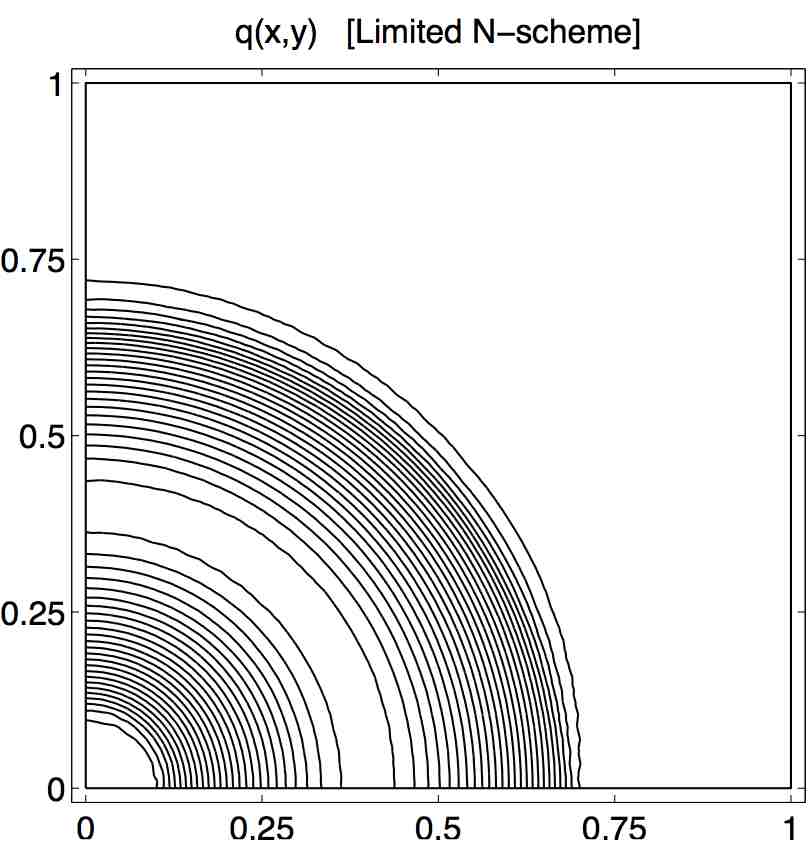} \,
(d)  \includegraphics[height=62mm]{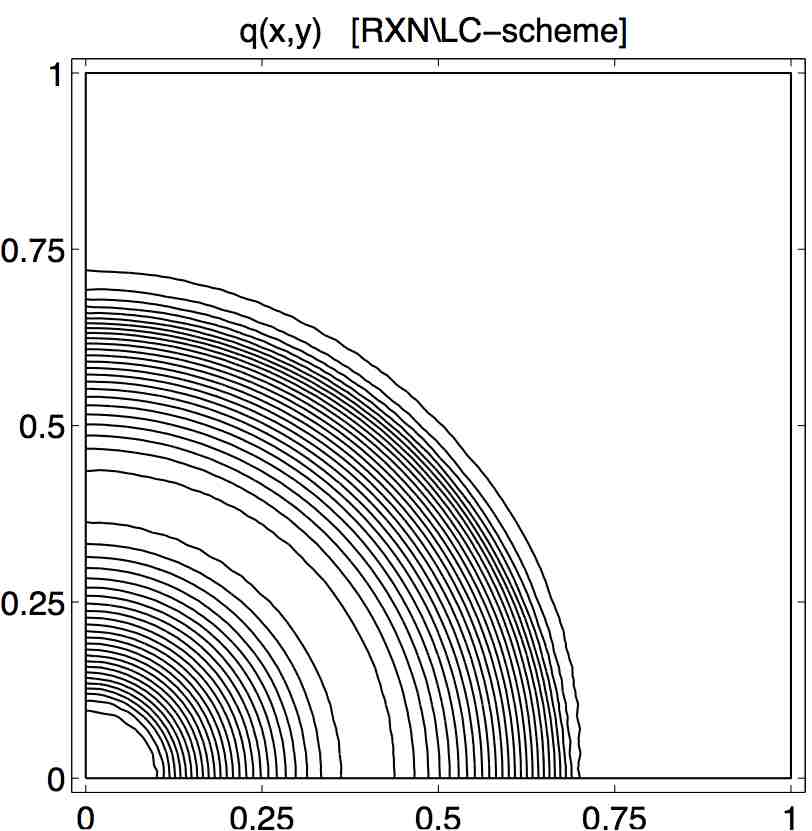}}

\end{center}
\caption{Advection equation example. Shown in these panels
are (a) the basic N-scheme, (b) the basic
RXN-scheme, (c)
the limited N-scheme (no convergence correction
is needed for the limited N-scheme on scalar equations), and
(d) the RXN$\backslash$LC-scheme 
(convergence corrections are needed for the limited RXN-scheme,
even for scalar problems). These results show that the 
basic RXN scheme is far more diffusive than the N-scheme.
However, with limiting and convergence corrections, the 
RXN$\backslash$LC gives results comparable to the limited
N-scheme.}
\label{fig:adv1}
\end{figure}

\begin{figure}[!t]
\begin{center}
\includegraphics[height=70mm]{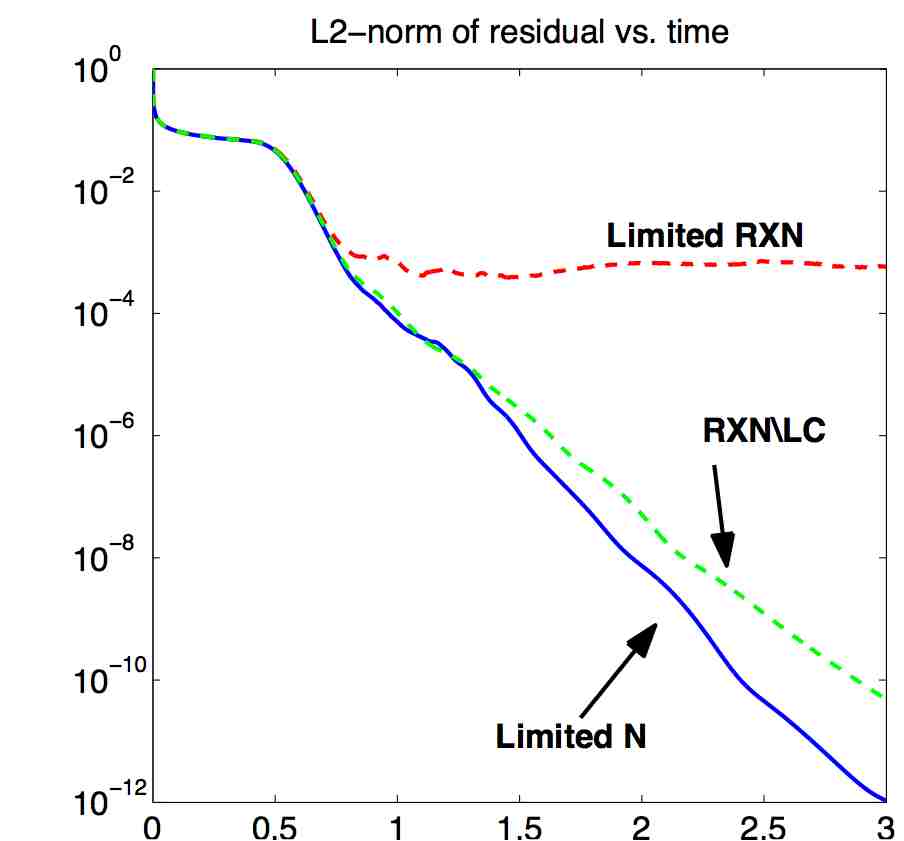}
\end{center}
\caption{Advection equation example.
$L_2$-norms of the total residual for the limited N-scheme,
the limited RXN-scheme, and the RXN$\backslash$LC-scheme. This
plot shows that the limited N-scheme does not need additional
convergence corrections for scalar equations, but the limited
RXN-scheme clearly does.}
\label{fig:adv2}
\end{figure}

Results on a grid with 5592 elements and
2903 nodes is shown in Figure \ref{fig:adv1}; displayed
in each panel are (a) the basic N-scheme, (b) the basic
RXN-scheme, (c) the limited N-scheme (no convergence correction
is needed for the limited N-scheme on scalar equations), and
(d) the RXN$\backslash$LC-scheme 
(convergence corrections are needed for the limited RXN-scheme,
even for scalar problems). These results show that the 
basic RXN scheme is far more diffusive than the N-scheme.
However, with limiting and convergence corrections, the 
RXN$\backslash$LC-scheme gives results comparable to the limited
N-scheme. Convergence histories for the 
limited N-scheme, limited RXN-scheme, and the RXN$\backslash$LC-scheme are shown in Figure \ref{fig:adv2}.

Clearly there is no advantage in using the RXN$\backslash$LC-scheme
over the limited N-scheme for a scalar problem, since the
two methods have the same computational cost and the scalar limited
N-scheme does not require convergence corrections. However,
for hyperbolic systems such as the Euler equations, the 
RXN$\backslash$LC-scheme provides a simpler algorithm
with lower computational cost than the N$\backslash$LC-scheme.

\subsection{Transonic flow around the NACA 0012 airfoil}
Next we consider transonic flow around the NACA 0012 airfoil
using the compressible Euler equations from gas dynamics
as our model. The Euler equations can be written as
\begin{equation}
  \begin{bmatrix}
     \rho \\ \rho \vec{u} \\ {\mathcal E}
  \end{bmatrix}_{,t} + \nabla \cdot 
   \begin{bmatrix}
     \rho \vec{u} \\ \rho \vec{u} \vec{u} + p \id \\
     \vec{u} \left( {\mathcal E} + p \right)
  \end{bmatrix} = 0,
\end{equation}
where $\rho$ is the fluid density, $\vec{u} = (u^1, \, u^2)$
is the fluid velocity, ${\mathcal E}$ is the total energy,
and $p$ is the fluid pressure. The ideal gas law closes
the system by relating the pressure to the other fluid variables:
\begin{equation}
  {\mathcal E} = \frac{p}{\gamma -1} + \frac{1}{2} \rho \| \vec{u} \|^2,
\end{equation}
where $\gamma$ is the ideal gas constant. In
this example we take $\gamma=1.4$.

\begin{figure}[!t]
\begin{center}
(a) \includegraphics[height=55mm]{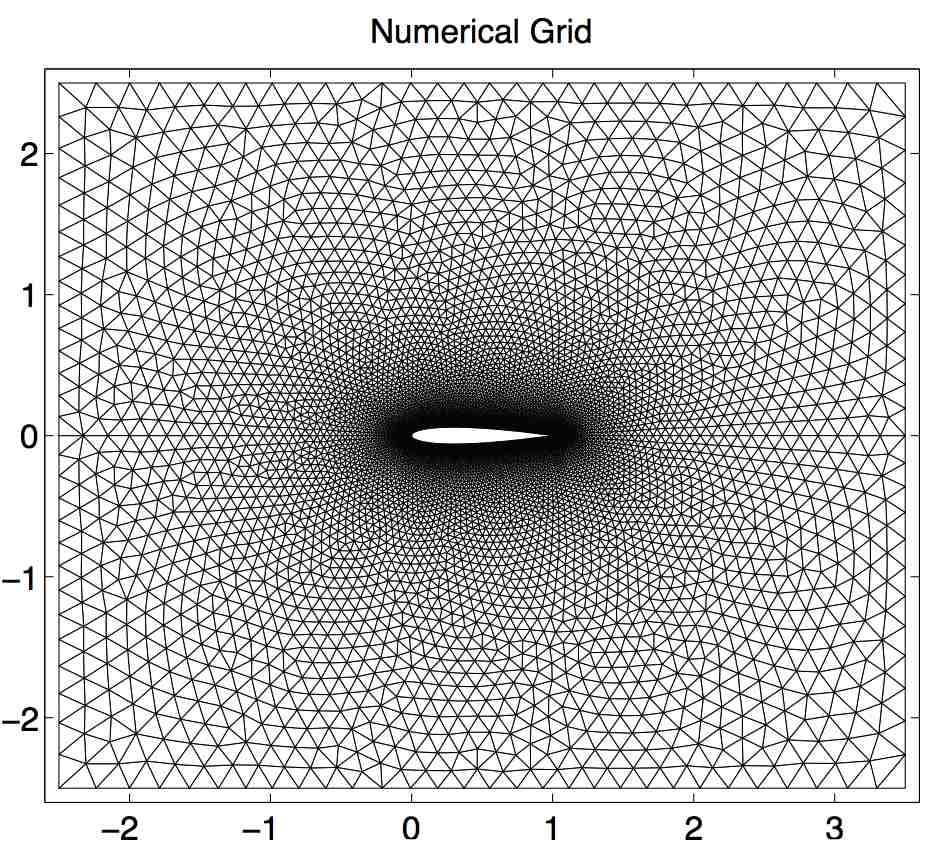}
\,
(b)\includegraphics[height=55mm]{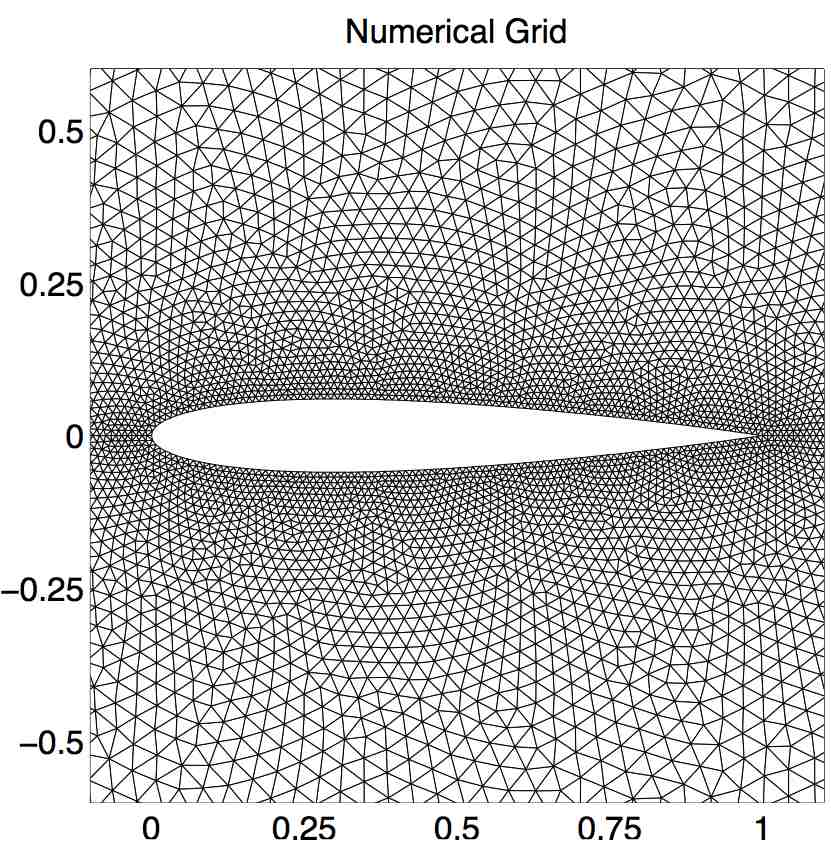}
\end{center}
\caption{The numerical grid for the NACA 0012 example. Panel
(a) shows the entire domain, while Panel (b) shows a zoomed-in
view of the airfoil. The grid has a total of 14,284 elements and
7,298 nodes. The smallest grid elements near the airfoil 
have a triangle radius $h$ that is roughly 30 times smaller than that of
the largest grid elements on the outer boundaries.}
\label{fig:naca_grid}
\end{figure}

The computational domain and numerical grid is shown
in Figure \ref{fig:naca_grid}(a). In Figure \ref{fig:naca_grid}(b), we show
a zoomed-in view of the numerical grid near the airfoil.
The boundary conditions are such that subsonic flow
with Mach number 0.85
enters from the left boundary at an angle of $+1^{\circ}$
from the horizontal axis. As the flow impinges on
the airfoil, two supersonic bubbles are created
above and below the airfoil. The supersonic flow
is decelerated to the ambient subsonic flow
through the creation of two shock waves (again, one
above the airfoil and one below).
This problem has been considered in several papers
including \cite{article:Ab01,article:Ab06}.

\begin{figure}[!t]
\begin{center}
(a) \includegraphics[height=80mm]{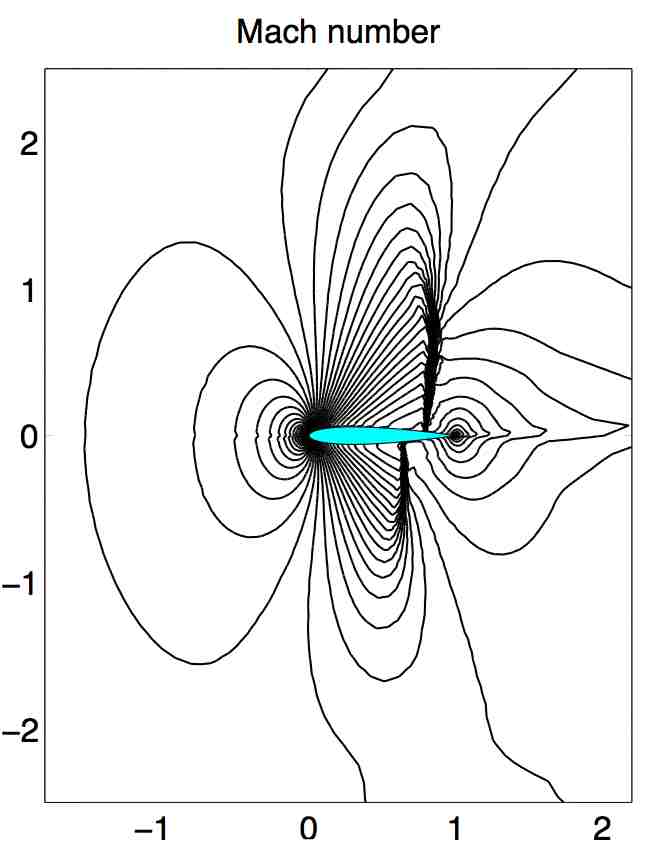}
\,
(b)\includegraphics[height=80mm]{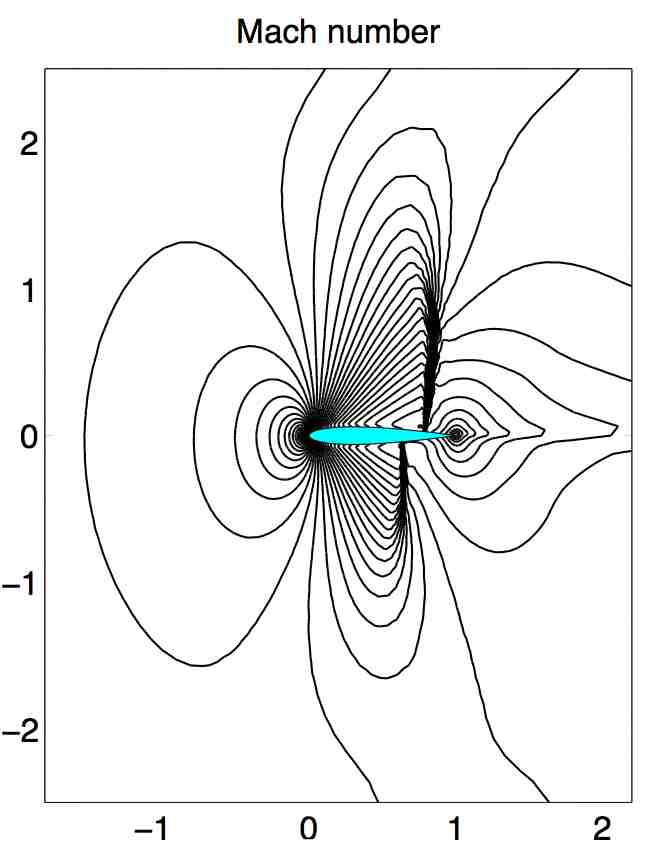}

\vspace{1mm}

(c) \includegraphics[height=80mm]{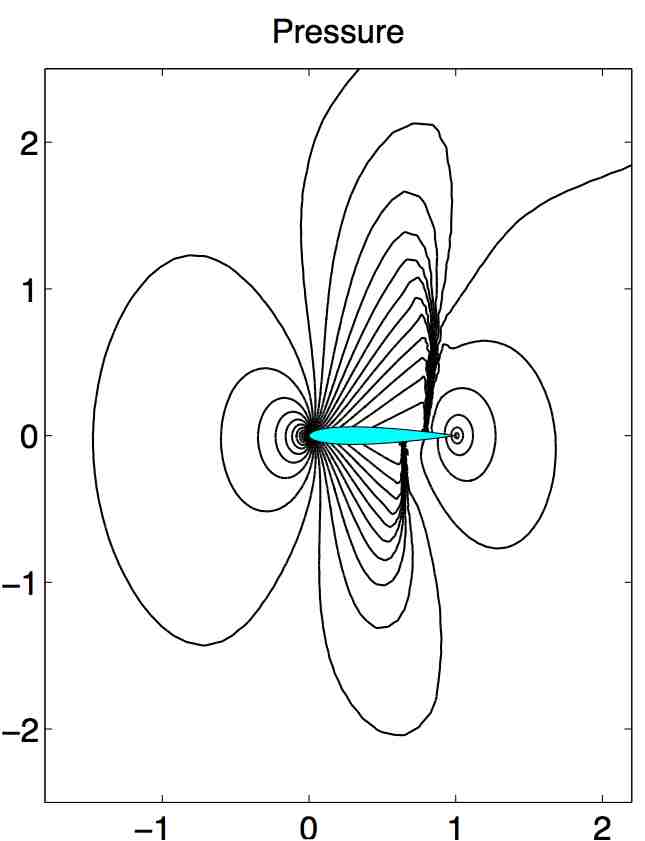}
\,
(d)\includegraphics[height=80mm]{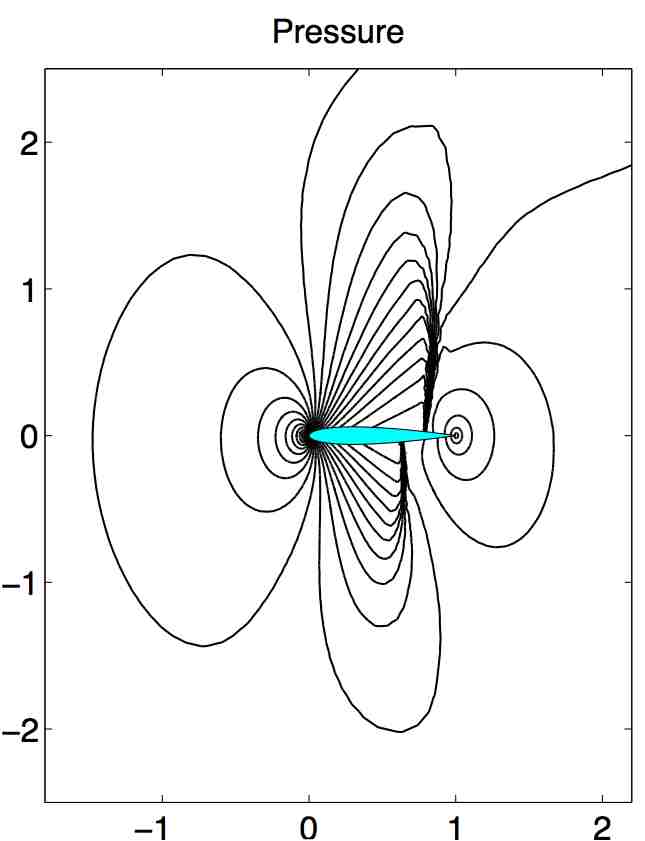}
\end{center}
\caption{Steady-state solution of transonic flow around the NACA 0012 airfoil.
Panels (a) and (b) show isolines of the Mach number for the
N$\backslash$LC-scheme and RXN$\backslash$LC-scheme, respectively.
Panels (c) and (d) show isolines of the pressure for the
N$\backslash$LC-scheme and RXN$\backslash$LC-scheme, respectively.
}
\label{fig:naca1}
\end{figure}

\begin{figure}[!t]
\begin{center}

(a) \includegraphics[angle=90,width=62mm]{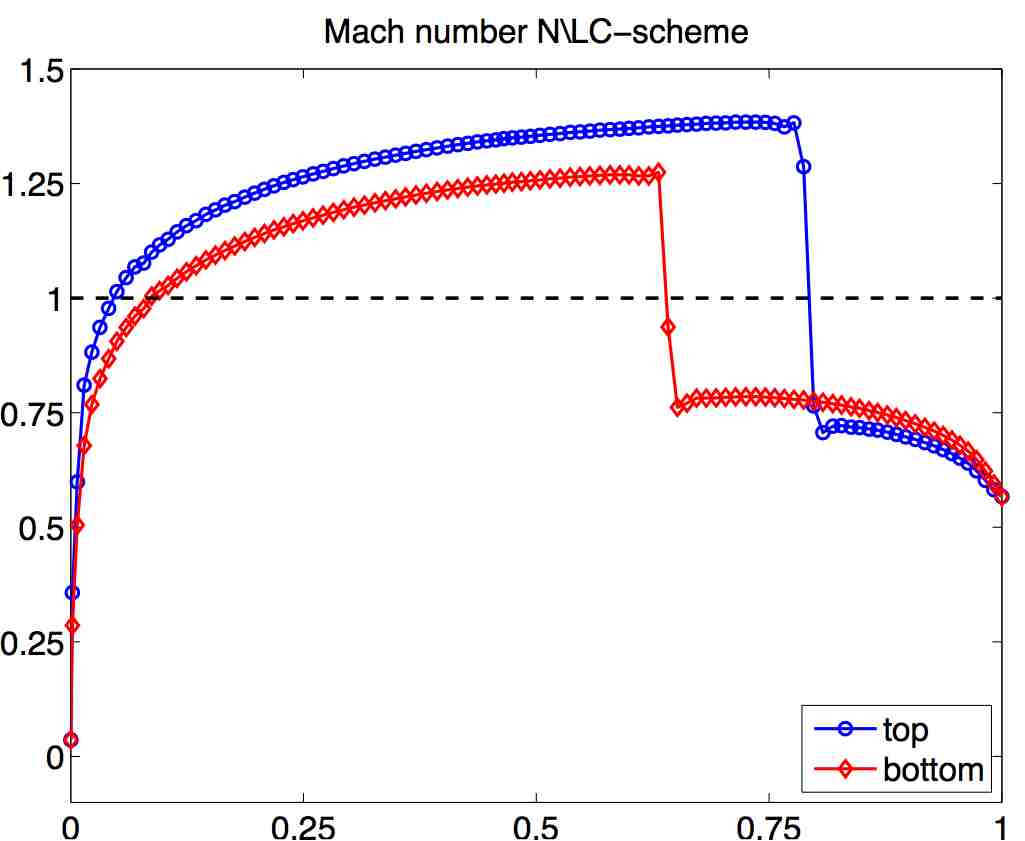}
\,
(b)\includegraphics[angle=90,width=62mm]{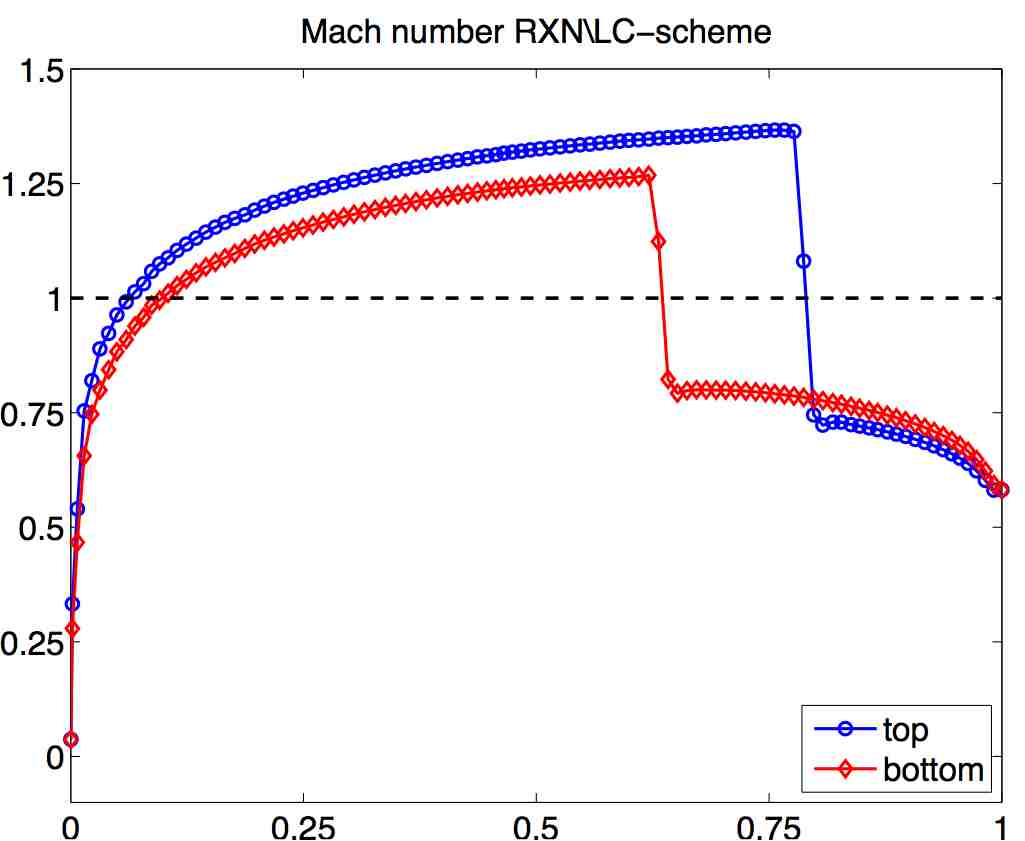}

\mbox{
(c) \includegraphics[angle=90,width=62mm]{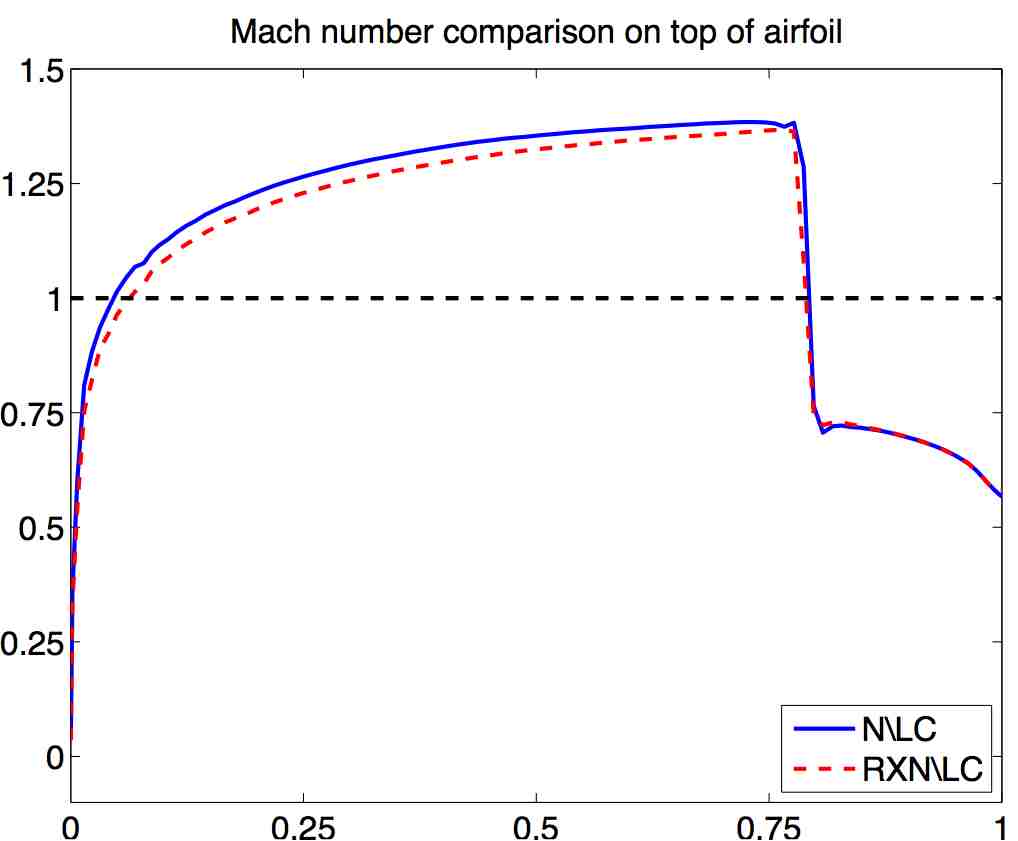}
\,
(d)\includegraphics[angle=90,width=62mm]{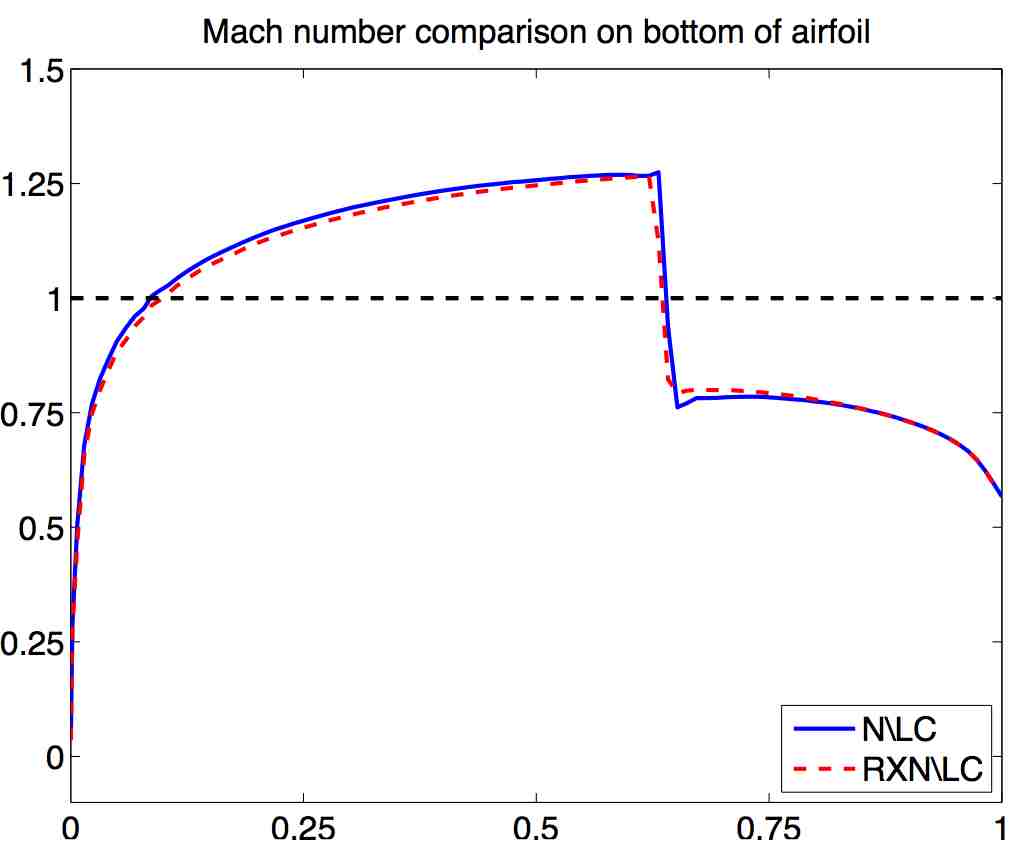} 
}

\end{center}
\caption{The Mach number along the top and bottom edges
of the NACA 0012 airfoil. Panel (a) is the N$\backslash$LC-scheme solution, while
Panel (b) is the RXN$\backslash$LC-scheme solution. 
In Panels (c) and (d) we directly compare the Mach number
profiles for each method: Panel (c) shows the Mach number
on the top portion of the airfoil and Panel (d) shows
the Mach number on the bottom portion of the airfoil.
These results show that the RXN$\backslash$LC-scheme
produces comparable results to the N$\backslash$LC-scheme.}
\label{fig:naca2}
\end{figure}

\begin{figure}[!t]
\begin{center}
(a)  \includegraphics[angle=90,height=85mm]
{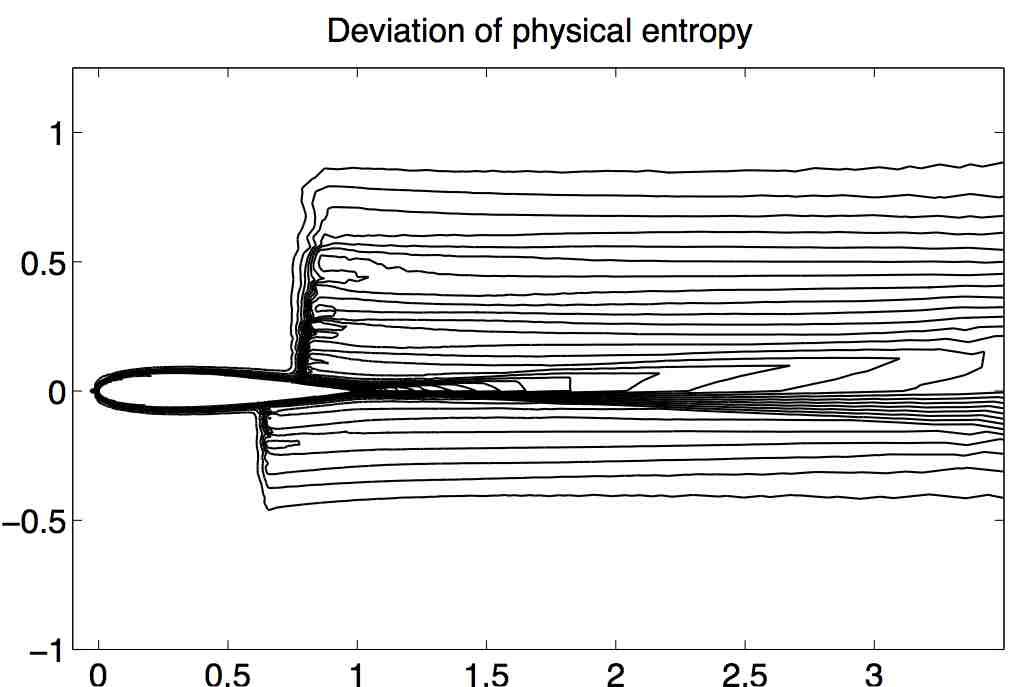} \, \, \,
(b)  \includegraphics[angle=90,height=85mm]{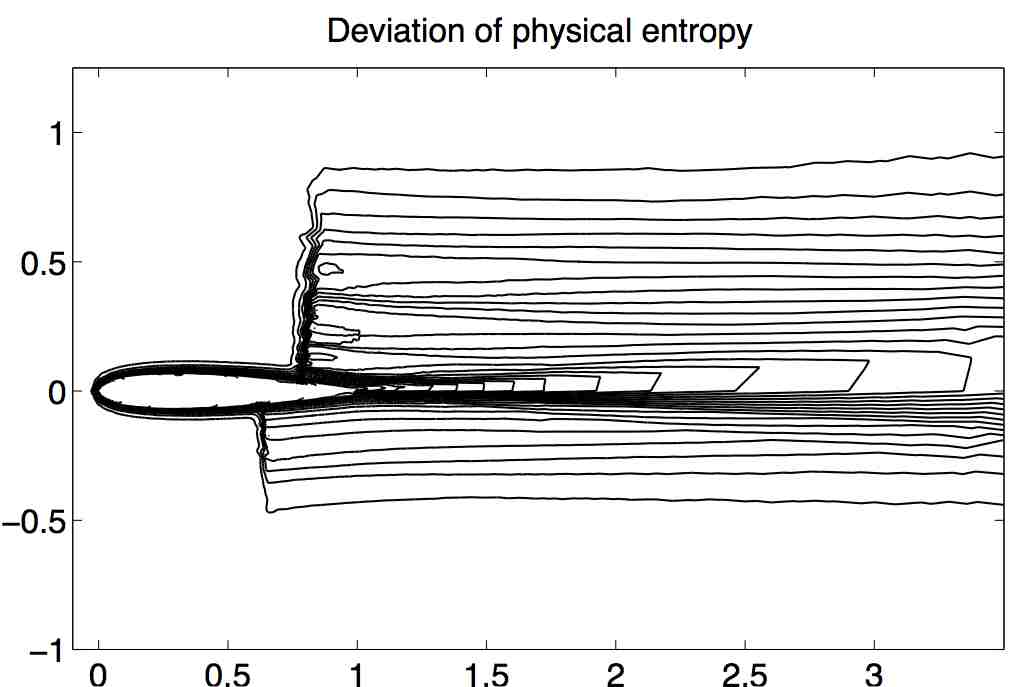}
\end{center}
\caption{Deviation of the physical entropy, 
$s = \log(p/\rho^{\gamma})$,  from the ambient
entropy, $s_{\infty} = \log(1/\gamma^\gamma)$:
$\Sigma = ( s - s_{\infty} )/|s_{\infty}|$. 
Panel (a) is the N$\backslash$LC-scheme and panel
(b) is the RXN$\backslash$LC-scheme. These
results again show that the RXN$\backslash$LC scheme
is slightly less accurate than the N$\backslash$LC scheme
near the airfoil. The same contour values are plotted
in each panel: $\Sigma = $ 0.002  :  0.002  :  0.08.
The minimum and maximum values of $\Sigma$
for the N$\backslash$LC and the RXN$\backslash$LC schemes
are $(-4.906 \times 10^{-4}, \, 8.029 \times 10^{-2})$ and
$(-2.778 \times 10^{-5}, \, 5.805 \times 10^{-2})$, respectively.}
\label{fig:naca_entropy}
\end{figure}

Isolines of the Mach number and the pressure are shown in Figure 
\ref{fig:naca1}.  We also plot in 
Figure \ref{fig:naca2} the Mach number along the
top and bottom edges of the airfoil. From this figure
we note that the location of the shocks is in very good agreement
between the two methods,
while the location at the front of the airfoil where the solution goes
from subsonic to supersonic is slightly different for each method.
Furthermore, the RXN$\backslash$LC-scheme is more diffusive
than the N$\backslash$LC-scheme, which results in slightly
more entropy production near the airfoil for RXN$\backslash$LC
than N$\backslash$LC. This can be seen both in Figure \ref{fig:naca1},
where we see bending of the Mach isolines near the airfoil, as
well as in Figure \ref{fig:naca_entropy}.
Overall, however, both of these figures indicate remarkable
agreement between the two solutions. In particular, the RXN$\backslash$LC
solution is far closer to the N$\backslash$LC solution than the
MUSCL-type scheme that was presented in \cite{article:Ab01}.
We also note that the RXN$\backslash$LC-scheme runs twice
as fast as the N$\backslash$LC-scheme.

\begin{figure}[!t]
\begin{center}
\includegraphics[height=60mm]{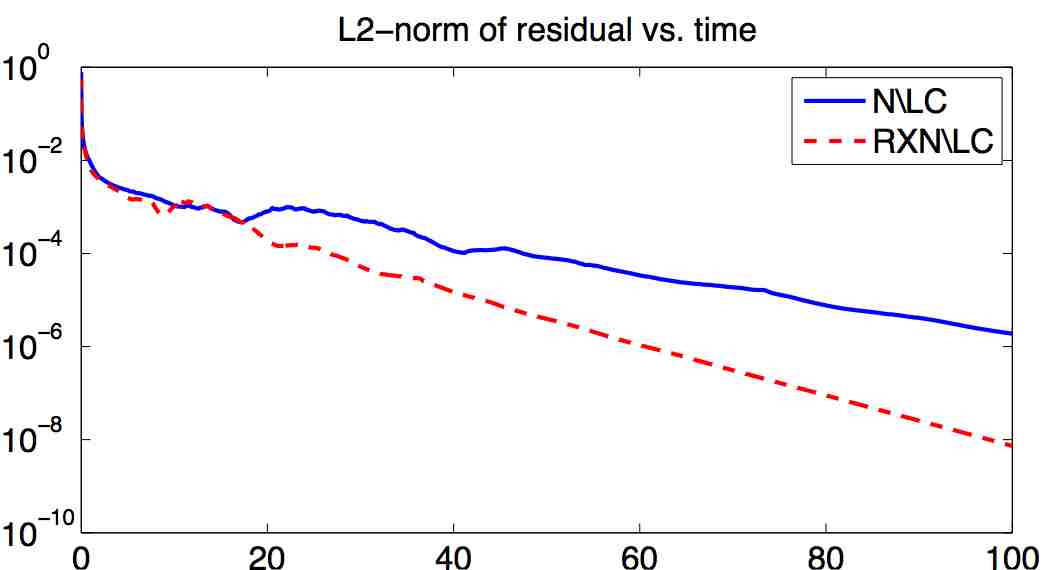}
\end{center}
\caption{$L_2$-norm of the total residual as a function of time for
the NACA 0012 problem.
	We note that the fix of Abgrall \cite{article:Ab06} (see
Section \ref{sec:Abfix}) is critically important
 	in bringing both methods to convergence. Without this fix
	both methods stall at a total residual of only about
	$10^{-2}$.}
\label{fig:naca_residual}
\end{figure}

Finally, in Figure \ref{fig:naca_residual} we show the
$L_2$-norm of the total residual as a function of time.
We note that the fix of Abgrall \cite{article:Ab06} (see
Section \ref{sec:Abfix}) is critically important
in bringing both methods to convergence. Without this fix
both methods stall at a total residual of only about
$10^{-2}$. We also find that the RXN$\backslash$LC-scheme
has a slightly better convergence rate than the N$\backslash$LC-scheme.

\subsection{Supersonic flow around a cylinder}
\begin{figure}[!t]
\begin{center}
(a) \includegraphics[width=35mm]{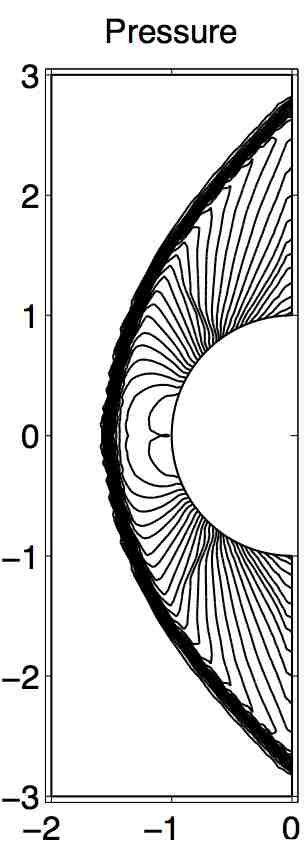}
\hspace{5mm}
(b) \includegraphics[width=35mm]{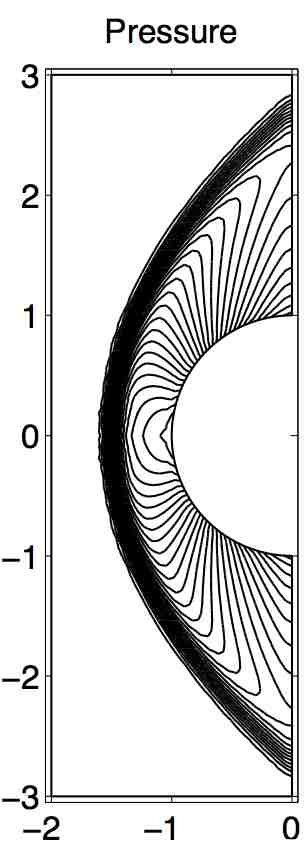}

(c) \includegraphics[width=35mm]{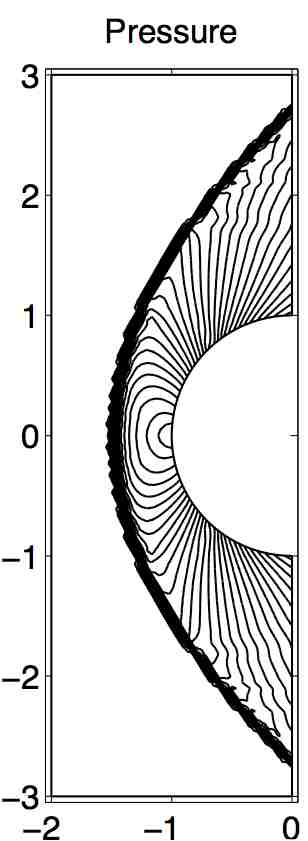}
\hspace{5mm}
(d) \includegraphics[width=35mm]{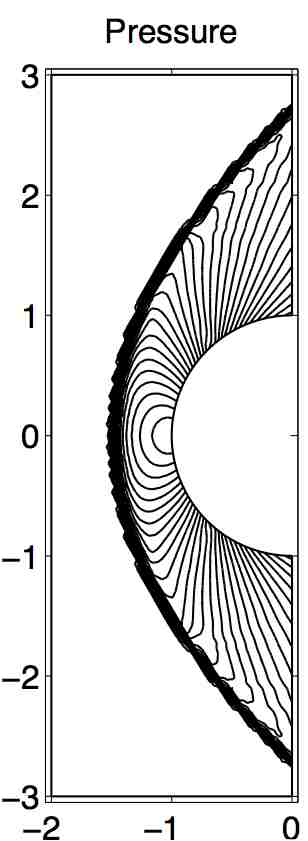}
\end{center}
\caption{Supersonic flow past a cylinder with ${\mathcal M}_{\infty} = 5$.
Shown in these panels are
the steady-state solutions as computed with the (a) basic N-scheme,
(b) basic RXN-scheme, (c) N$\backslash$LC-scheme, and
(d) RXN$\backslash$LC-scheme. These results show that the
RXN solution is much more diffusive than the N-scheme solution;
however, once the limiters and convergence corrections are
included, the N$\backslash$LC and  RXN$\backslash$LC schemes
produce comparable results.}
\label{fig:supersonic}
\end{figure}

\begin{figure}[!t]
\begin{center}
\includegraphics[height=60mm]{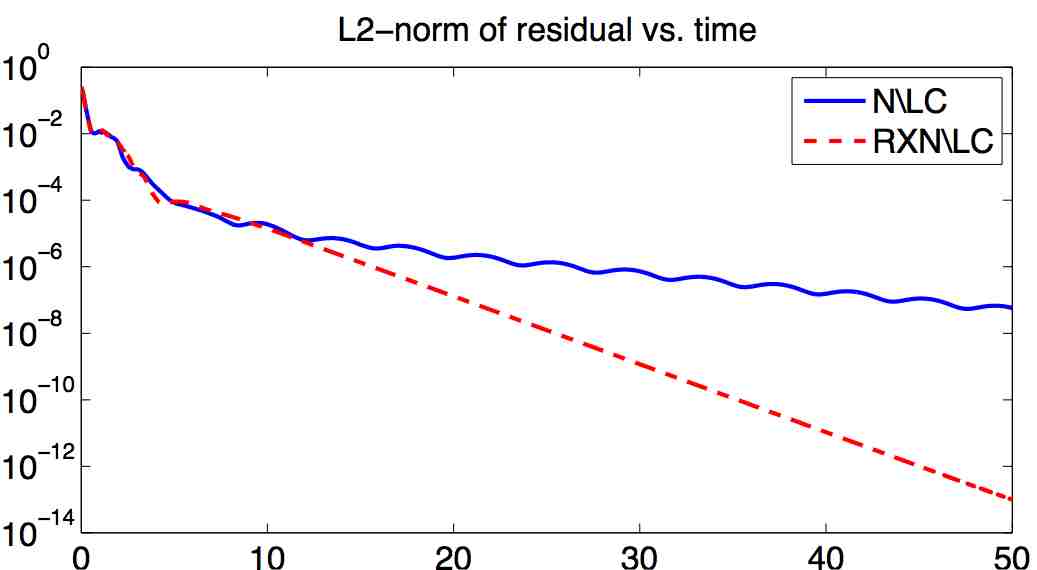}
\end{center}
\caption{$L_2$-norm of the total residual as a function of time for
supersonic flow past a cylinder.}
\label{fig:supersonic_residual}
\end{figure}

Next we consider flow around a cylinder with Mach number
${\mathcal M}_{\infty} = 5$. The computational domain is
$[-2, 0] \times [-3, 3]$ with a cylinder of unit radius centered at $(0,0)$.
In this example we found that we needed to run all the schemes 
at a CFL number of 0.4 in order to obtain well-converged results.
The steady-state pressure on a grid with 5144 elements and 2656 nodes
is shown in Figure \ref{fig:supersonic} for the following schemes:
(a) the  basic N-scheme, (b) the basic RXN-scheme, (c) the
N$\backslash$LC-scheme, and (d) the RXN$\backslash$LC-scheme.
The basic RXN-scheme solution is far more diffusive than the basic N-scheme,
but once limiters and the convergence corrections are added, the
N$\backslash$LC and RXN$\backslash$LC schemes produce
comparable results. In fact, the RXN$\backslash$LC-scheme
converges faster than the N$\backslash$LC-scheme, as can
be seen in the convergence plot in Figure \ref{fig:supersonic_residual}.
Finally, we note that the RXN$\backslash$LC-scheme again runs
about twice as fast as the N$\backslash$LC-scheme.

\subsection{Subsonic flow around a cylinder}
\begin{figure}[!t]
\begin{center}

(a) \includegraphics[width=95mm]{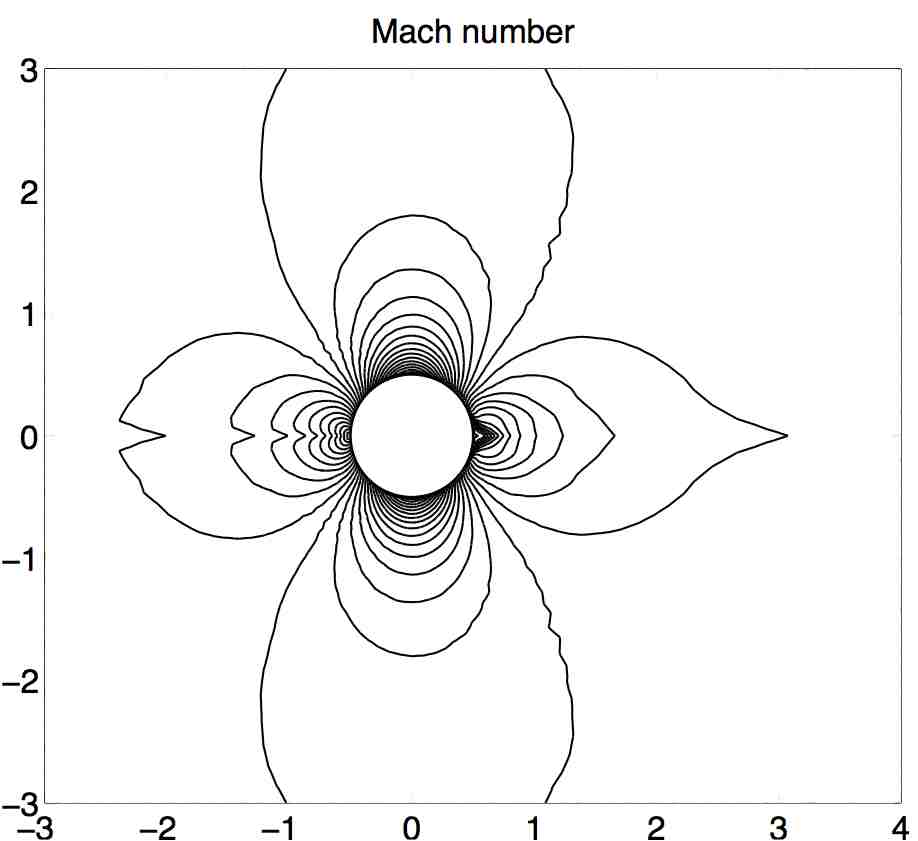}

(b) \includegraphics[width=95mm]{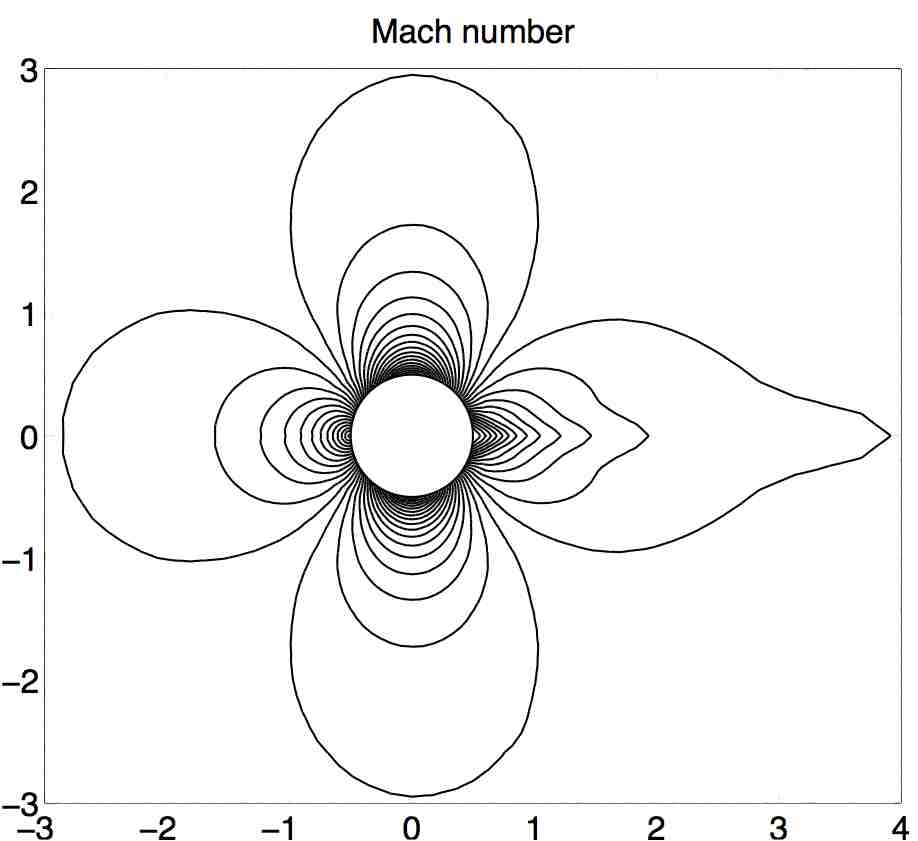}
\end{center}
\caption{Subsonic flow past a cylinder with ${\mathcal M}_{\infty} = 0.35$.
Shown in these panels are isolines of the Mach number for
the (a) N$\backslash$LC and (b) RXN$\backslash$LC schemes. 
Near the cylinder both methods produce comparable results. Away from the cylinder the grid resolution becomes coarser; and therefore, visible differences
in the two methods appear. In these regions the
RXN$\backslash$LC scheme produces
slightly more diffused contours than the N$\backslash$LC scheme.}
\label{fig:subsonic_mach}
\end{figure}

Finally, we consider flow around a cylinder with Mach number
${\mathcal M}_{\infty} = 0.35$. 
This problem has been considered in several papers
including \cite{article:Ab01,article:Ab06}.
The computational domain is
$[-7,7] \times [-7, 7]$ with a cylinder of radius $1/2$ centered at $(0,0)$.
The steady-state Mach number on a grid with 12552 elements and 6404 nodes
is shown in Figure \ref{fig:subsonic_mach} for the (a) 
N$\backslash$LC and (b) RXN$\backslash$LC schemes.
Near the cylinder both methods produce comparable results. Away from the cylinder the grid resolution becomes coarser; and therefore, visible differences
in the two methods appear. In these regions the
RXN$\backslash$LC scheme produces
slightly more diffused contours than the N$\backslash$LC scheme.

\begin{figure}[!t]
\begin{center}
(a) \includegraphics[width=130mm]{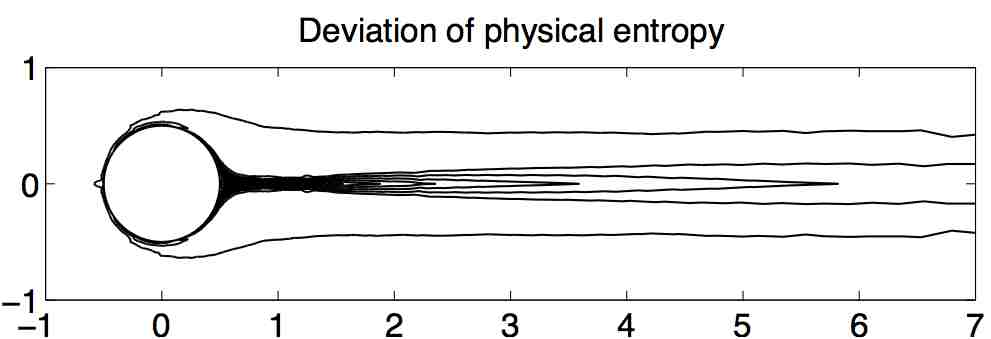}

(b) \includegraphics[width=130mm]{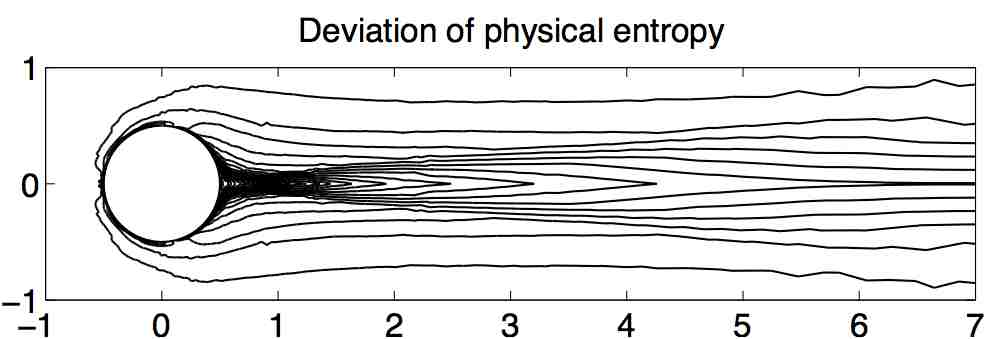}
\end{center}
\caption{Subsonic flow past a cylinder with ${\mathcal M}_{\infty} = 0.35$.
Shown in these panels are  the deviation of the physical entropy, 
$s = \log(p/\rho^{\gamma})$,  from the ambient
entropy, $s_{\infty} = \log(1/\gamma^\gamma)$:
$\Sigma = ( s - s_{\infty} )/|s_{\infty}|$. 
Panel (a) is the N$\backslash$LC-scheme and panel
(b) is the RXN$\backslash$LC-scheme. 
The minimum and maximum values of $\Sigma$
for the N$\backslash$LC and the RXN$\backslash$LC schemes
are $(-4.101 \times 10^{-3}, \, 4.324 \times 10^{-2})$ and
$(-1.471 \times 10^{-3}, \, 1.291 \times 10^{-2})$, respectively.
Each panel consists of 31 contour lines ranging from the minimum
 to the maximum $\Sigma$ for each scheme. Therefore, these results show that
the RXN$\backslash$LC-scheme has a smaller entropy deviations, but that this error is more spread out behind the cylinder, while the
N$\backslash$LC-scheme has
larger entropy deviations, but that this error is more concentrated
near the $x$-axis.}
\label{fig:subsonic_entropy}
\end{figure}

Shown in Figure \ref{fig:subsonic_entropy} 
are  the deviation of the physical entropy, 
$s = \log(p/\rho^{\gamma})$,  from the ambient
entropy, $s_{\infty} = \log(1/\gamma^\gamma)$:
$\Sigma = ( s - s_{\infty} )/|s_{\infty}|$. 
Panel (a) is the N$\backslash$LC-scheme and panel
(b) is the RXN$\backslash$LC-scheme. 
The minimum and maximum values of $\Sigma$
for the N$\backslash$LC and the RXN$\backslash$LC schemes
are $(-4.101 \times 10^{-3}, \, 4.324 \times 10^{-2})$ and
$(-1.471 \times 10^{-3}, \, 1.291 \times 10^{-2})$, respectively.
Each panel consists of 31 contour lines ranging from the minimum
 to the maximum $\Sigma$ for each scheme. Therefore, these results show that
the RXN$\backslash$LC-scheme has a smaller entropy deviations, but that this error is more spread out behind the cylinder, while the
N$\backslash$LC-scheme has
larger entropy deviations, but that this error is more concentrated
near the $x$-axis.

\begin{figure}[!t]
\begin{center}
\includegraphics[height=60mm]{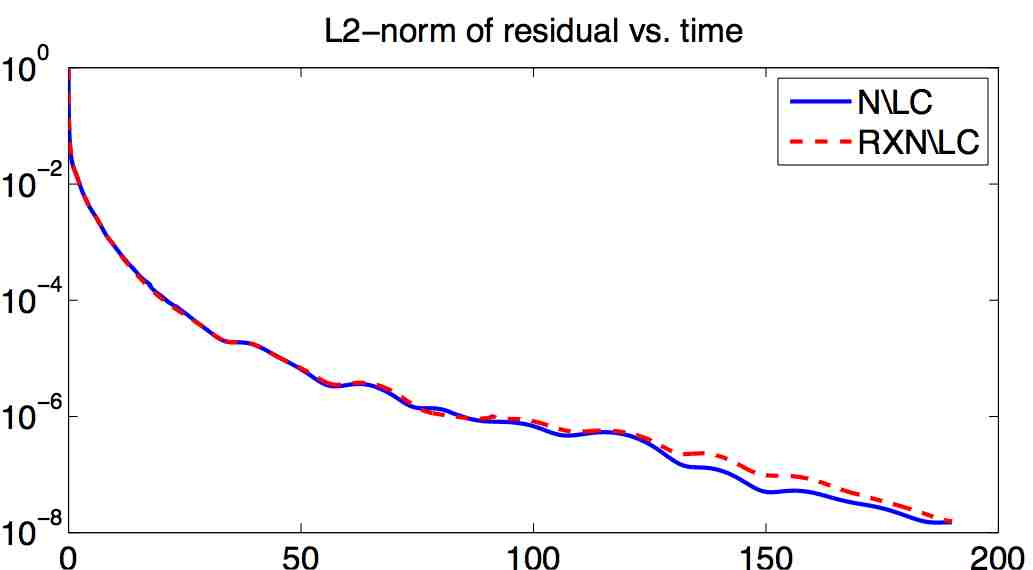}
\end{center}
\caption{$L_2$-norm of the total residual as a function of time for
subsonic flow past a cylinder. Both methods give essentially the
same convergence rates for this example.}
\label{fig:subsonic_residual}
\end{figure}

The total residual as a function of time is shown in Figure
\ref{fig:subsonic_residual}. Both methods give essentially the
same convergence rates for this example.
Finally, we note that the RXN$\backslash$LC-scheme again runs
about twice as fast as the N$\backslash$LC-scheme.

\section{Conclusions}
In this work we have extended the results of LeVeque
and Pelanti \cite{article:LePe01} to genuinely multidimensional
residual distribution schemes. Specifically, we have shown 
that the N-scheme, both the scalar and the systems version,
can be derived from a relaxation principle. Furthermore, 
using a genuinely multidimensional extension of
the 1D local Lax-Friedrichs relaxation principle, we have
derived a novel residual distribution scheme. The main
benefit of this approach is that it does not require
the inversion of an $m\times m$ matrix, where $m$ is
the number of conserved variables, at each time
step in each grid element. The new method also does
not require the use of Roe-Struijs-Deconinck averages.
Using several examples of the 2D Euler equations from
gas dynamics, including an
example of transonic flow around the NACA 0012 airfoil, 
supersonic flow around a cylinder, and subsonic flow
around a cylinder, we have compared the limited and corrected
N-scheme (N$\backslash$LC) with the newly proposed
scheme (RXN$\backslash$LC). These comparisons show that
despite being computationally less expensive, the new method is
capable of producing results comparable to those
of the N$\backslash$LC-scheme, although often with slightly
more numerical diffusion. For more complicated equations such
as magnetohydrodynamics or the general relativistic Einstein equations,
we believe that the benefit of a simpler and computationally less expensive
algorithm will far outweigh the slight increase in numerical dissipation.
We will consider some of these more complicated systems in
future work.

Finally, we would like to point out that the numerical code used in this work, including all of the numerical grids, 
will be made publicly
available as part of the {\tt REDPACK} software project.
For more information see
\[
\text{{\tt http://www.math.wisc.edu/$\sim$rossmani/software.html}} \, .
\]

{\bf Acknowledgements. } 
The author would like to thank the two anonymous reviewers
for their very helpful comments. This
work was supported in part by NSF
grants DMS--0619037
and DMS--0711885.

\bibliographystyle{plain}

\end{document}